\documentclass[12pt]{article}
\usepackage{epsfig}
\usepackage{color}
\usepackage[T2A]{fontenc}
\usepackage[utf8]{inputenc}
\usepackage[english,russian]{babel}
\usepackage{hyperref}
\usepackage{amssymb,amsfonts,amsmath,mathtext,cite,enumerate,float} 
\usepackage{graphicx}
\usepackage{algorithm}
\usepackage{algorithmic}
\usepackage{indentfirst}
\usepackage{amsthm}
\usepackage{epigraph}
\usepackage{xcolor}
\usepackage{verbatim}
\usepackage{url}

\begin{document}

\begin{center}
\large \bf О построении градиентного метода квадратичной оптимизации, оптимального с точки зрения минимизации расстояния до точного решения 
\end{center}
\vspace{18pt}

\large
\bf \copyright 2025 г. \phantom{20} Н.В. Плетнев \vspace{12pt}

\normalsize
\rm
\begin{center}\it
141701, г. Долгопрудный, Институтский пер., д. 9, МФТИ, Россия

e–mail: nikita\_pletnev@list.ru

\end{center}

\footnotesize
 \vspace{6pt}
Задачи квадратичной оптимизации в гильбертовом пространстве часто возникают при решении некорректных задач для дифференциальных уравнений. При этом известно целевое значение функционала. Кроме того, структура функционала позволяет вычислять градиент с помощью решения корректных задач, что позволяет применять методы первого порядка. Данная статья посвящена построению $m$-моментного метода минимальных ошибок --- эффективного метода, минимизирующего расстояние до точного решения. Доказывается сходимость и оптимальность построенного метода, а также невозможность равномерной сходимости методов, работающих в подпространствах Крылова. Проводятся численные эксперименты, демонстрирующие эффективность применения $m$-моментного метода минимальных ошибок к решению различных некорректных задач: начально-краевой задачи для уравнения Гельмгольца, ретроспективной задачи Коши для уравнения теплопроводности, обратной задачи термоакустики.
 \vspace{6pt}

 {\bf Ключевые слова:} некорректные и обратные задачи, квадратичная оптимизация, оптимизация в гильбертовом пространстве, минимизация расстояния до точного решения, начально-краевая задача для уравнения Гельмгольца, ретроспективная задача Коши для уравнения теплопроводности, обратная задача термоакустики.
\normalsize

\vspace{12pt}
\centerline{1. Введение} 
\vspace{12pt}
Теория некорректных и обратных задач --- очень важный раздел математической физики. Область их применения в физике, геофизике, технике, медицине весьма обширна. Многие обратные задачи могут быть представлены в виде операторного уравнения:
\begin{equation}\label{Aq=f}
    Aq=f
\end{equation}

Здесь оператор $A$ определён в некотором гильбертовом пространстве $H$, решение $q^*$ предполагается единственным, а $f$ --- результат измерений, доступных в условиях задачи. При этом вычисление результата действия оператора $A$ на произвольный элемент $H$ --- корректная, прямая задача. 

Примеры задач, представимых в таком виде: начально-краевая задача для уравнения Гельмгольца, задача Коши для уравнения теплопроводности, интегральное уравнение Фредгольма 1-ого рода.

Часто оператор $A$ не является линейным из-за неоднородных граничных условий. При рассмотрении таких задач вводится оператор $A_0$, который определяется аналогично, но с нулями вместо всех начальных и граничных условий, кроме заданного функцией $q$. Важное свойство: $Aq-Aq'=A_0(q-q').$

Широкое распространение получил подход к решению подобных задач с помощью сведения их к задачам выпуклой оптимизации: 
\begin{equation}\label{Jq}
    J(q)=\frac{1}{2}||Aq-f||^2\rightarrow\min
\end{equation} 

Его эффективность обусловлена сравнительной простотой вычисления градиента: 
\begin{equation}
    \nabla J(q) = A^*(Aq-f)
    \label{nabla_Jq}
\end{equation}
Здесь $A^*$ --- сопряжённый оператор для $A_0$.

Градиент удовлетворяет условию Липшица с константой $L=||A_0||^2$. Поскольку задача вычисления оператора $A$ является корректной, $A_0$ ограничен. 
\begin{equation}
    ||\nabla J(q) - \nabla J(q')||=||A^*(Aq-f) - A^*(Aq'-f)||=||A^*A_0(q-q')||\leq L||q-q'||.
\end{equation}

Часто для решения используется метод градиентного спуска, также применимы ускоренные методы, в том числе и метод сопряжённых градиентов. Структура функционала позволяет выполнить точную вспомогательную одномерную минимизацию для произвольного направления. В некоторых случаях возможна также минимизация расстояния до точного решения. Построению эффективного метода оптимизации, основанного на этом подходе, и посвящена статья.

\vspace{12pt}
\centerline{2. Работы, основанные на похожих идеях} 
\vspace{12pt}

В статье \cite{polyak1969step} был предложен, а в статье \cite{devanathan2024pmm} развит развит так называемый PMM --- метод минорант Поляка. Он основан на следующей идее: на каждом шаге $k$ подбирается <<миноранта>> $m_k(q)$ --- функция простого вида, для которой $m_k(q)\leq J(q)$ и $m_l(q_k)=J(q_k)$. Подбирается точка $q_{k+1}$, для которой $m_k(q_{k+1})=J(q^*)=0$. При этом $J(q_{k+1})<J(q_k)$.

Поскольку <<естественная>> миноранта для гладкого выпуклого функционала --- это его линейное приближение, метод принимает вид:
\begin{equation}
    \label{polyak_method}
    q_{k+1} = q_k - \frac{J(q_k)}{||\nabla J(q_k)||^2} \nabla J(q_k)
\end{equation}

В дальнейшем будет показано, что одномерная минимизация расстояния до точного решения приводит к похожей формуле.

Эксперименты показывают, что метод Поляка эффективнее одномерной минимизации расстояния. Как указывает сам автор, его скорость сходимости примерно соответствует характеристикам обычных градиентных методов. Однако, оказывается возможным достижение лучших результатов при использовании многомерной вспомогательной минимизации.

В статье \cite{goujaud2022quadratic} предпринята попытка построить метод квадратичной минимизации (адаптивный тяжёлый шарик), оптимальный в смысле минимизации расстояния до точного решения. Однако для такого метода это расстояние должно монотонно убывать. Эксперименты показывают, что построенный в указанной статье метод не обладает этим свойством. Однако данная работа решает проблему: метод строится, а его свойства доказываются.

\vspace{12pt}
\centerline{3. Построение оптимального метода квадратичной оптимизации} 
\vspace{12pt}
\centerline{\it 3.1. Одномерная минимизация расстояния до точного решения: метод минимальных ошибок}
\vspace{12pt}

Специальный вид функционала и градиента позволяет (согласно \cite{kabanikhin2012book}, раздел 2.6) выбрать длину шага в общем случае. Пусть шаг $k$ делается в направлении $s_k$: 
\begin{equation}q_{k+1}=q_k+\alpha_k s_k.\label{step}
\end{equation}

$$\alpha_k(\rho)=\arg\min\limits_{\alpha\geq 0} ||q_{k} + \alpha s_k - q^*||^2.$$

Элементарные преобразования минимизируемого выражения позволяют получить искомое значение коэффициента:\\
\begin{center}
$||q_{k} + \alpha s_k - q^*||^2 = ||q_{k} - q^*||^2 + 2\alpha \langle q_k-q^*,\ s_k\rangle+\alpha^2||s_k||^2\Rightarrow$ 
\begin{equation}\Rightarrow \alpha_k(\rho) = -\frac{\langle q_k-q^*,\ s_k\rangle}{||s_k||^2}.
\label{alpha(rho)}
\end{equation}
\end{center}

Поскольку точное решение $q^*$ неизвестно, эта формула применима лишь тогда, когда скалярное произведение в числителе может быть преобразовано к вычислимому виду. Соответственно, минимизация расстояния до точного решения возможна не при любом выборе направления спуска.

Если она возможна, то на шаге $k$ достигается уменьшение расстояния до точки минимума:
$$||q_{k+1}-q^*||^2 = ||q_{k} - q^*||^2 - \frac{\langle q_k-q^*,\ s_k\rangle^2}{||s_k||^2} = ||q_{k} - q^*||^2 - ||q_{k+1} - q_{k}||^2.$$

Суммируя от $0$ до $n$, получаем: $$||q_{0} - q^*||^2-||q_{n+1} - q^*||^2 = \sum\limits_{k=0}^{n}||q_{k+1} - q_{k}||^2.$$
Убывающая и ограниченная снизу нулём последовательность $||q_{n} - q^*||^2$ имеет предел, перейдём к нему:
\begin{equation}||q_{0} - q^*||^2 - \lim\limits_{n\rightarrow\infty} ||q_{n} - q^*||^2 = \sum\limits_{k=0}^{\infty}||q_{k+1} - q_{k}||^2.
\label{limNorm2}
\end{equation}

Следовательно:
\begin{itemize}
    \item ряд в правой части всегда сходится;
    \item его сумма не превосходит $||q_{0} - q^*||^2$;
    \item равенство суммы $||q_{0} - q^*||^2$ --- необходимое и достаточное условие сильной сходимости метода с шагами $s_k$ и длиной каждого шага $\alpha_k(\rho)$;
    \item $\lim\limits_{k\rightarrow\infty}||q_{k+1}-q_k|| = 0$ (необходимое условие сходимости ряда). Из признака сравнения рядов следует, что для любой числовой последовательности $x_k$ такой, что ряд $\sum\limits_{k=1}^{\infty}x_k$ расходится, выполнено условие $||q_{k+1}-q_k||=o(\sqrt{x_k})$. В частности, $||q_{k+1}-q_k||=o\left(\frac{1}{\sqrt{k}}\right)$.
\end{itemize} 

Если $s_k=-\nabla J(q_k)$, то $$\alpha_k(\rho) = -\frac{\langle q_k-q^*,\ s_k\rangle}{||s_k||^2} = \frac{\langle q_k-q^*,\ \nabla J(q_k)\rangle}{||\nabla J(q_k)||^2} = \frac{\langle q_k-q^*,\ A^*(Aq_k-f)\rangle}{||\nabla J(q_k)||^2} =$$ $$= \frac{\langle A_0(q_k-q^*),\ Aq_k-f\rangle}{||\nabla J(q_k)||^2} = \frac{||Aq-f||^2}{||\nabla J(q_k)||^2}\Rightarrow$$ 
\begin{equation}\Rightarrow\alpha_k(\rho)=\frac{2J(q_k)}{||\nabla J(q_k)||^2}.
\label{alpha(rho)grad}
\end{equation}

\vspace{12pt}
\centerline{\it 3.2. Моментный метод минимальных ошибок}
\vspace{12pt}

Моментными называют методы следующего вида:

\begin{equation}
    q_{k+1}=q_k - \alpha_k\nabla J(q_k) +\gamma_k(q_k-q_{k-1})
    \label{moment_dq}
\end{equation}
или
\begin{equation}
    s_k=-\nabla J(q_k) + \beta_k s_{k-1};\quad q_{k+1}=q_k + \alpha_ks_k
    \label{moment_s}
\end{equation}.

Поскольку векторы $s_{k-1}$ и $q_k-q_{k-1}$ коллинеарны, эти записи эквивалентны, и одни коэффициенты можно выразить через другие: $\gamma_k\alpha_{k-1}=\alpha_k\beta_k$.

Можно заметить, что к данному классу относятся различные уже рассмотренные методы, такие как метод сопряжённых градиентов, метод подобных треугольников или метод тяжёлого шарика. Идея заключается в том, чтобы оптимально подобрать коэффициенты.

Выполним вспомогательную двумерную минимизацию расстояния до точного решения при вычислении $q_{k+1}$ в форме (\ref{moment_dq}). Пусть $\gamma_0=0$, $\alpha_0=\alpha_0(\rho)$. $q_k\not= q^*$ (иначе новый шаг не нужен).

$$||q_{k+1}-q^*||^2 = ||q_k-q^* - \alpha_k\nabla J(q_k) +\gamma_k(q_k-q_{k-1})||^2=$$
$$=||q_k-q^*||^2+\alpha_k^2||\nabla J(q_k)||^2+\gamma_k^2||q_k-q_{k-1}||^2-2\alpha_k\langle q_k-q^*, \nabla J(q_k)\rangle +$$$$+ 2\gamma_k\langle q_k-q^*, q_k-q_{k-1}\rangle-2\alpha_k\gamma_k\langle \nabla J(q_k), q_k-q_{k-1}\rangle.$$

Вычислим частные производные. Их равенство нулю является необходимым условием экстремума.
$$
    \frac{1}{2}\frac{\partial}{\partial\alpha_k} ||q_{k+1}-q^*||^2 = \alpha_k||\nabla J(q_k)||^2 -\gamma_k\langle \nabla J(q_k), q_k-q_{k-1}\rangle - \langle q_k-q^*, \nabla J(q_k)\rangle;
$$
$$
    \frac{1}{2}\frac{\partial}{\partial\gamma_k} ||q_{k+1}-q^*||^2 = -\alpha_k\langle \nabla J(q_k), q_k-q_{k-1}\rangle + \gamma_k||q_k-q_{k-1}||^2 + \langle q_k-q^*, q_k-q_{k-1}\rangle.
$$

Преобразуем слагаемые, содержащие $q^*$:
$$\langle q_k-q^*, \nabla J(q_k)\rangle = \langle q_k-q^*, A^*(Aq_k-f)\rangle = \langle A_0(q_k-q^*), Aq_k-f\rangle =$$
\begin{equation}
    = ||Aq_k-f||^2 = 2J(q_k);
    \label{sqal_q_qs_nabla}
\end{equation}
$$\langle q_k-q^*, q_k-q_{k-1}\rangle = 0,$$
потому что $q_k$ выбрано на предыдущем шаге, как ближайший к $q^*$ элемент линейного многообразия $q_{k-1}+Lin\left\{ \nabla J(q_{k-1}), q_{k-1}-q_{k-2}\right\}$, а это проекция на гиперплоскость.

Итак, получаем систему линейных уравнений:
$$
\begin{cases}
    \alpha_k||\nabla J(q_k)||^2 -\gamma_k\langle \nabla J(q_k), q_k-q_{k-1}\rangle = 2J(q_k);\\
    -\alpha_k\langle \nabla J(q_k), q_k-q_{k-1}\rangle + \gamma_k||q_k-q_{k-1}||^2 = 0.
\end{cases}$$

Для её решения необходимо проверить, что она не является вырожденной. Для этого вычислим определитель:
$$\Delta = ||\nabla J(q_k)||^2||q_k-q_{k-1}||^2 - \langle \nabla J(q_k), q_k-q_{k-1}\rangle^2.$$

Неравенство Коши --- Буняковского показывает, что всегда $\Delta\geq 0,$ и равенство эквивалентно коллинеарности $\nabla J(q_k)$ и $q_k-q_{k-1}$. Поскольку $\langle q_k-q^*, q_k-q_{k-1}\rangle = 0,$ одновременно выполняется и следующее равенство: $$2J(q_k) = \langle q_k-q^*, \nabla J(q_k)\rangle = \langle q_k-q^*, \lambda(q_k-q_{k-1})\rangle =0.$$ Но это возможно только в том случае, когда $q_k$ --- точное решение задачи. По исходному предположению, это не так. 

Близость $\Delta$ к нулю обозначает плохую обусловленность системы и низкую точность очередного шага в условиях реальных вычислений. Этот признак может быть использован для формулировки критерия остановки метода.

Если решение задачи оптимизации не достигнуто, то $\Delta>0$. Тогда система имеет единственное решение $(\alpha_k, \gamma_k)$. Оно является точкой минимума $||q_{k+1}-q^*||^2$, поскольку матрица вторых производных --- положительно определённая по критерию Сильвестра: её главные миноры $||\nabla J(q_k)||^2$ и  $\Delta$ положительны.

Осталось вычислить это решение. Из второго уравнения: $$\gamma_k = \alpha_k\cdot \frac{\langle \nabla J(q_k), q_k-q_{k-1}\rangle}{||q_k-q_{k-1}||^2}.$$
Подставляя в первое уравнение, получаем:
$$\alpha_k = \frac{2J(q_k)}{||\nabla J(q_k)||^2 - \frac{\langle \nabla J(q_k), q_k-q_{k-1}\rangle^2}{||q_k-q_{k-1}||^2}}.$$

Полученный метод является оптимальным среди моментных методов по признаку уменьшения расстояния до точного решения после одного шага.

Введём обозначение: $$s_k= - \nabla J(q_k) + \frac{\langle \nabla J(q_k), q_k-q_{k-1}\rangle}{||q_k-q_{k-1}||^2}(q_k - q_{k-1}).$$
Поскольку $q_k - q_{k-1} = \alpha_{k-1}s_{k-1},$ выполняется равенство: 
$$s_k = -\nabla J(q_k) + \beta_ks_{k-1},$$
где $$\beta_k = \frac{\langle \nabla J(q_k), s_{k-1}\rangle}{||s_{k-1}||^2}.$$

Очевидно, $\langle s_k, q_k-q_{k-1}\rangle=0$. Поэтому $$||s_k||^2 + \frac{\langle \nabla J(q_k), q_k-q_{k-1}\rangle^2}{||q_k-q_{k-1}||^2} = ||\nabla J(q_k)||^2,$$ откуда $$\alpha_k=\frac{2J(q_k)}{||s_k||^2}.$$

Вычисление $\alpha_k(\rho)$ для введённого шага $s_k$ приводит к тому же результату:
$$\alpha_k(\rho) = -\frac{\langle q_k-q^*, s_k\rangle}{||s_k||^2} = \frac{\langle q_k-q^*, \nabla J(q_k)\rangle}{||s_k||^2} - \beta_k\frac{\langle q_k-q^*, s_{k-1}\rangle}{||s_k||^2} = \frac{2J(q_k)}{||s_k||^2},$$
так как $\langle q_k-q^*, s_{k-1}\rangle=0$ по доказанному в ходе преобразования частных производных.

Итак, построенный локально оптимальный моментный метод может быть представлен как моментный метод (\ref{moment_s})
с коэффициентом сопряжённости 

\begin{equation}
    \beta_k = \frac{\langle \nabla J(q_k), s_{k-1}\rangle}{||s_{k-1}||^2},
    \label{beta_k_2}
\end{equation} вычисленным из соображения ортогональности соседних шагов, и выбором длины шага 
\begin{equation}
    \alpha_k=\alpha_k(\rho) = \frac{2J(q_k)}{||s_k||^2}
    \label{alpha_k_2}
\end{equation} по принципу минимизации расстояния до точного решения.

\vspace{12pt}
\centerline{\it 3.3. Многошаговый метод минимальных ошибок}
\vspace{12pt}

Рассмотрим следующее обобщение построенного метода, зависящее от натурального параметра $m$. Назовём его $m$-моментным методом минимальных ошибок. При $m=1$ получится метод из раздела 3.2.

\begin{equation}
q_{k+1} = \arg\min\limits_{q\in q_k + Lin\{-\nabla J(q_k), q_k-q_{k-1},\ldots,q_{k-m+1}-q_{k-m}\}} ||q-q^*||^2.
    \label{m_moment_dq}
\end{equation}

То есть, на каждом шаге выбирается ближайшая к решению точка из линейного многообразия, проходящего через последнюю полученную точку и натянутого на векторы антиградиента и $m$ предыдущих шагов.

Если $k<m$, то используются все $k$ шагов, в этом случае размерность пространства равна $k+1$.

Для краткости обозначим $h_j:=q_{j+1}-q_j$ при всех $j$.
\begin{equation}
    q_{k+1} = q_k -\gamma_k^0\nabla J(q_k) + \sum\limits_{i=1}^m\gamma_k^i h_{k-i}.
    \label{m_moment_h}
\end{equation}
Коэффициенты $\gamma_k^j$ в (\ref{m_moment_h}) получаются путём решения задачи минимизации (\ref{m_moment_dq}).

Выпишем минимизируемое выражение:
$$||q_{k+1}-q^*||^2 = ||q_k-q^*||^2 + \left(\gamma_k^0\right)^2||\nabla J(q_k)||^2 + \sum\limits_{i=1}^m \left(\gamma_k^i\right)^2||h_{k-i}||^2 -$$$$- 2\gamma_k^0\langle q_k-q^*, \nabla J(q_k)\rangle  + 2\sum\limits_{1\leq i<j\leq m}\gamma_k^i\gamma_k^j\langle h_{k-i}, h_{k-j}\rangle -$$$$- 2\sum\limits_{i=1}^m\gamma_k^0\gamma_k^i\langle\nabla J(q_k),h_{k-i}\rangle + 2\sum\limits_{i=1}^m\gamma_k^i\langle q_k-q^*, h_{k-i}\rangle.$$

Преобразуем слагаемые, содержащие неизвестное $q^*$:
\begin{equation}
    \langle q_k-q^*, \nabla J(q_k)\rangle = \langle A_0(q_k-q^*), Aq_k-f\rangle = ||Aq_k-f||^2=2J(q_k);\label{sqal_q_qs_nabla_multi}
\end{equation}
\begin{equation}
    \langle q_k-q^*, h_{k-i}\rangle = \langle q_k-q_{k-i+1}, h_{k-i}\rangle + \langle q_{k-i+1}-q^*, h_{k-i}\rangle.
    \label{sqal_q_qs_h}
\end{equation}

Второе слагаемое в (\ref{sqal_q_qs_h}) равно $0$, поскольку точки $q^*$, $q_{k-i}$, $q_{k-i+1}$ образуют прямоугольный треугольник с прямым углом $q_{k-i+1}$. Это связано с тем, что выбранная на шаге $k-i$ точка является ближайшей к $q^*$ в своём подпространстве --- значит и на принадлежащей ему прямой.

Теперь (\ref{sqal_q_qs_h}) можно преобразовать:
\begin{equation}
    \langle q_k-q^*, h_{k-i}\rangle = \langle q_k-q_{k-i+1}, h_{k-i}\rangle = \langle \sum\limits_{j=1}^{i-1} h_{k-j}, h_{k-i}\rangle = \sum\limits_{j=1}^{i-1}\langle h_{k-j}, h_{k-i}\rangle.
    \label{sqal_dq_h}
\end{equation}

Все слагаемые, кроме $||q_k-q^*||^2$, не зависящего от вектора $\gamma_k$, удалось выразить в вычислимом виде. Следовательно, минимум может быть найден с помощью необходимого (для выпуклой квадратичной функции --- достаточного) условия экстремума. Для этого вычислим частные производные по $\gamma_k^i$:

$$\frac{1}{2}\frac{\partial}{\partial \gamma_k^0} ||q_{k+1}-q^*||^2 = \gamma_k^0||\nabla J(q_k)||^2 - \sum\limits_{i=1}^m\gamma_k^i\langle\nabla J(q_k),h_{k-i}\rangle - 2J(q_k);$$
при $1\leq i\leq m$
$$\frac{1}{2}\frac{\partial}{\partial \gamma_k^i} ||q_{k+1}-q^*||^2 = \gamma_k^i||h_{k-i}||^2 + \sum\limits_{j=1,\ j\not=i}^m\gamma_k^j\langle h_{k-j},h_{k-i}\rangle -$$$$- \gamma_k^0\langle\nabla J(q_k),h_{k-i} + \sum\limits_{j=1}^{i-1}\langle h_{k-j}, h_{k-i}\rangle.$$

Таким образом, условие экстремума можно представить в матричном виде: $$G\gamma_k=b,$$ где $$G_{00} = ||\nabla J(q_k)||^2,\ G_{0i}=G_{i0}=-\langle \nabla J(q_k), h_{k-i}\rangle,$$$$G_{ii} = ||h_{k-i}||^2,\ G_{ij}=G_{ji} = \langle h_{k-i},h_{k-j}\rangle,$$
$$\gamma_k=\left(\gamma_k^0,\gamma_k^1,\ldots,\gamma_k^m\right)^T,$$
$$b=\left(2J(q_k), 0, -\langle h_{k-1}, h_{k-2}\rangle, \ldots,-\sum\limits_{j=1}^{i-1}\langle h_{k-i}, h_{k-j}\rangle, \ldots,-\sum\limits_{j=1}^{m-1}\langle h_{k-m}, h_{k-j}\rangle\right)^T.$$

Это представление вызывает следующие вопросы. Может ли матрица $G$ оказаться вырожденной? Возможно ли упрощение вычислений (решения системы линейных уравнений порядка $m+1$)? Ответ на эти вопросы даёт следующая лемма.

{\bf Лемма 1.}
\label{lemma:ortogonality}
Каждый шаг ортогонален $m$ предыдущим (или всем, если его номер меньше $m$), то есть $\langle h_k, h_{k-i}\rangle = 0$ при $1\leq i\leq m$.

{\bf Доказательство.}
Применим метод математической индукции по номеру шага $k$. База очевидна: при $k=0$ или $k=1$ и произвольном $m\geq 1$ метод совпадает с моментным методом из раздела 3.2, для которого ортогональность соседних шагов уже была доказана.

Пусть лемма верна для всех шагов с номерами меньше $k$. Тогда на шаге $k$ все компоненты вектора $b$, кроме имеющей нулевой номер, равны нулю по предположению индукции. Также равны нулю все элементы матрицы $G$, кроме нулевой строки, нулевого столбца и главной диагонали.

Поскольку $h_k = -\gamma_k^0\nabla J(q_k) + \sum\limits_{i=1}^m\gamma_k^i h_{k-i}$, при  $1\leq i\leq m$. можно вычислить скалярное произведение:
$$\langle h_k, h_{k-i}\rangle = -\gamma_k^0\langle\nabla J(q_k), h_{k-i}\rangle + \gamma_k^i ||h_{k-i}||^2 = \sum\limits_{j=0}^m G_{ji}\gamma_k^j = b_i=0.$$

QED. $ $

Из леммы 1 и формулы (\ref{sqal_dq_h}) сразу следует, что при $1\leq i\leq m$ 
$$\langle q_k-q^*, h_{k-i}\rangle = 0.$$
Предположим теперь, что система векторов, из линейной оболочки которых выбирается вектор очередного шага, является линейно зависимой, то есть $$\nabla J(q_k) = \sum\limits_{i=1}^m\lambda_ih_{k-i}.$$ Умножим скалярно на $q_k-q^*$, получим равенство $$2J(q_k)=\sum\limits_{i=1}^m\lambda_i\langle q_k-q^*, h_{k-i}\rangle =0,$$ которое выполняется только при достижении точного решения.

Соответственно, при $q_k\not=q^*$ эта система векторов линейно независима. Для решения системы уравнений $G\gamma_k=b$ выпишем и преобразуем (методом Гаусса) её расширенную матрицу:
$$\left\|\left.
\begin{array}{cccc}
    ||\nabla J(q_k)||^2 & -\langle\nabla J(q_k), h_{k-1}\rangle & \cdots & -\langle\nabla J(q_k), h_{k-m}\rangle\\
    -\langle\nabla J(q_k), h_{k-1}\rangle & ||h_{k-1}||^2 & \cdots & 0\\
    -\langle\nabla J(q_k), h_{k-2}\rangle & 0 & \cdots & 0\\
    \vdots & \vdots & \ddots & \vdots\\
    -\langle\nabla J(q_k), h_{k-m}\rangle & 0 & \cdots & ||h_{k-m}||^2
\end{array}
\right|
\begin{array}{c}
   2J(q_k) \\
    0 \\ 0 \\ \vdots \\ 0 
\end{array}\right\|\sim$$
$$\sim
\left\|\left.
\begin{array}{cccc}
    ||\nabla J(q_k)||^2 - \sum\limits_{i=1}^m\frac{\langle\nabla J(q_k), h_{k-i}\rangle^2}{||h_{k-i}||^2} & 0 & \cdots & 0\\
    -\langle\nabla J(q_k), h_{k-1}\rangle & ||h_{k-1}||^2 & \cdots & 0\\
    -\langle\nabla J(q_k), h_{k-2}\rangle & 0 & \cdots & 0\\
    \vdots & \vdots & \ddots & \vdots\\
    -\langle\nabla J(q_k), h_{k-m}\rangle & 0 & \cdots & ||h_{k-m}||^2
\end{array}
\right|
\begin{array}{c}
   2J(q_k) \\
    0 \\ 0 \\ \vdots \\ 0 
\end{array}\right\|$$

Можно заметить, что вычитаемая сумма квадратов в левом верхнем углу --- это квадрат модуля проекции $-\nabla J(q_k)$ на $Lin\{h_{k-1},\ldots,h_{k-m}\}$. Обозначим угол между антиградиентом и этой проекцией --- $\varphi_k$. Тогда левый верхний элемент матрицы равен $||\nabla J(q_k)||^2\sin^2{\varphi_k}$, и можно вычислить компоненты вектора $\gamma_k$:
$$\gamma_k^0=\frac{2J(q_k)}{||\nabla J(q_k)||^2\sin^2{\varphi_k}};\quad \gamma_k^i = \gamma_k^0\cdot\frac{\langle\nabla J(q_k), h_{k-i}\rangle}{||h_{k-i}||^2} \mbox{ при } 1\leq i\leq m.$$

Теперь шаг метода может быть представлен в следующем виде:
\begin{equation}
    s_k = -\nabla J(q_k) + \sum\limits_{i=1}^m \frac{\langle\nabla J(q_k), h_{k-i}\rangle}{||h_{k-i}||^2}h_{k-i}
    \label{s_m_moment}
\end{equation}
--- проекция антиградиента на ортогональное дополнение линейной оболочки $m$ предыдущих шагов; длина шага $\alpha_k$ вычисляется по формуле (\ref{alpha_k_2}).

Формула для $s_k$ может быть переписана рекуррентно:
\begin{equation}
    s_k = -\nabla J(q_k) + \sum\limits_{i=1}^m \frac{\langle\nabla J(q_k), s_{k-i}\rangle}{||s_{k-i}||^2}s_{k-i}.
    \label{s_m_moment_requrrent}
\end{equation}

При $k<m$ суммирование выполняется до индекса $i=k$. Формально положив $m=\infty$, получаем метод, каждый шаг которого ортогонален всем предыдущим. Его применение с $n$ шагами требует $n$ вычислений градиента и порядка $n^2$ вычислений скалярного произведения. Важно заметить, что, в отличие от вычисления градиента, скалярное произведение --- <<дешёвая>> операция, поэтому на практике квадратичная сложность не приводит к заметному увеличению времени работы.

\vspace{12pt}
\centerline{\it 3.4. Применение при составном функционале}
\vspace{12pt}

Все вычисления проведены в предположении, что функционал имеет простой квадратичный вид (\ref{Jq}). Но в некоторых задачах функционал более сложен: он состоит из нескольких квадратичных слагаемых с общей точкой минимума. Это, например, обратная задача термоакустики или задача восстановления источников для уравнения теплопроводности.

Такой функционал имеет следующий вид:
\begin{equation}
    J(q) = \frac{1}{2}\sum\limits_{l=1}^s ||A_l q - f_l||^2,
    \label{Jq_complex}
\end{equation}
причём $\forall l\in \{1,\ldots,m\}\; A_l q^* = f_l,$ поэтому $J(q^*)=0.$

Его градиент вычисляется аналогично простому случаю:
\begin{equation}
    \nabla J(q) = \sum\limits_{l=1}^s A_l^*(A_l q - f_l).
\end{equation}

В этом случае все преобразования этого раздела повторяются дословно, за исключением (\ref{sqal_q_qs_nabla}) и (\ref{sqal_q_qs_nabla_multi}):
$$\langle q_k-q^*, \nabla J(q_k)\rangle = \left\langle q_k-q^*, \sum\limits_{l=1}^s A_l^*(A_lq_k-f_l)\right\rangle =$$$$= \sum\limits_{l=1}^s\left\langle A_{l0}(q_k-q^*), A_lq_k-f_l\right\rangle = \sum\limits_{l=1}^s||A_lq_k-f_l||^2 = 2J(q_k).$$

Результат аналогичный, поэтому все полученные формулы сохраняются. Оценки сходимости также будут общими.

\vspace{12pt}
\centerline{\it 4. Теоремы о сходимости построенных методов}
\vspace{12pt}

{\bf Теорема 1. Об оптимальности $\infty$-моментного ММО.}

    $\infty$-моментный метод минимальных ошибок является глобально оптимальным для решения задачи минимизации квадратичного функционала в гильбертовом пространстве среди методов, работающих в подпространствах Крылова. То есть, никакой метод из указанного класса не позволяет достичь меньшего расстояния до точного решения при одинаковом количестве шагов.

{\bf Доказательство.}
    Обозначим $B=A^*A_0$ --- самосопряжённый оператор, связанный с задачей. Пусть $q_0$ --- фиксированный элемент множества, на котором определён оператор $A$, (если $A=A_0$, то $q_0=0$). Тогда градиент функционала представляется так: $$\nabla J(q) = B(q-q^*) = B(q-q_0)+A^*(Aq_0-f)=B(q-q_0)+\nabla J(q_0).$$

    Функционал можно представить в виде, типичном для применения классического метода сопряжённых градиентов: $$J(q)=\frac{1}{2}\left\langle q-q_0, B(q-q_0)\right\rangle - \left\langle A^*(f-Aq_0), q-q_0\right\rangle + ||Aq_0-f||^2.$$

    То есть, обозначая градиент в начальной точке $\nabla J(q_0)=g,$ получаем подпространства Крылова для задачи: $$K_n=Lin\{g, Bg,\ldots, B^{n-1}g\}.$$

    Докажем по индукции следующее утверждение о последовательности $q_n$, порождаемой $\infty$-моментным минимальных ошибок: $$\text{если векторы } g, Bg, B^2g,\ldots, B^ng \text{ линейно независимы, то  } \nabla J(q_n)\in K_{n+1}\backslash K_n, h_n\in K_{n+1}\backslash K_n.$$

    При $n=0$ оно очевидно: по определению $K_0=\{0\}$, $K_1=Lin\{\nabla J(q_0)\}$ --- одномерное подпространство ($g\not=0$, поскольку система, содержащая нулевой вектор, была бы линейно зависимой). 

    Пусть для всех шагов с номерами меньше $n$ утверждение выполнено. Докажем его для номера $n$. 
    $$\nabla J(q_n) = \nabla J(q_{n-1}) + B(q_n - q_{n-1})=\nabla J(q_{n-1}) + Bh_{n-1}.$$
    
    Первое слагаемое, по предположению индукции, лежит в $K_n\subset K_{n+1}$. Посмотрим на второе слагаемое. По предположению индукции $h_{n-1}\in K_n$ и $h_{n-1}\not\in K_{n-1}$, то есть $h_{n-1}$ раскладывается по системе $\{g,Bg,\ldots,B^{n-1}g\}$, причём последний коэффициент ненулевой, поэтому $Bh_{n-1}$ раскладывается по системе $\{Bg,B^2g,\ldots,B^ng\}$, и последний коэффициент --- снова ненулевой (разложение по линейно независимой системе единственно), то есть $Bh_{n-1}\in K_{n+1}$ и $Bh_{n-1}\not\in K_n$. Следовательно, $\nabla J(q_n)\in K_{n+1}$ (как сумма двух его элементов) и $\nabla J(q_n)\not\in K_n$ (иначе $Bh_{n-1}$ лежало бы в $K_n$, как разность его элементов).
    $$h_n = -\nabla J(q_n) + \sum\limits_{i=1}^n \gamma_n^i h_{n-i}.$$
    
    Здесь все слагаемые принадлежат $K_{n+1}$ (по доказанному), а все, кроме первого, --- $K_n$ (по предположению индукции). Отсюда сразу получается требуемое.

    Соответственно, на каждом шаге очередная точка $q_{n+1}$ выбирается из $q_n+K_{n+1}$. Поскольку $q_n-q_0\in K_{n+1}$, это множество совпадает с $q_0+K_{n+1}$.

    В рассматриваемом множестве $\infty$-моментный метод выбирает точку, ближайшую к решению. Поэтому совпадение множеств $$q_n + Lin\{-\nabla J(q_n), h_{n-1},\ldots,h_1, h_0\}=q_0 + K_{n+1}$$ и обозначает справедливость утверждения теоремы в случае, если для любого номера шага $n$ векторы $g, Bg, B^2g,\ldots, B^ng$ образуют линейно независимую систему.

    Если же для некоторого $n$ эта система впервые оказалась линейно зависимой, то максимальная размерность подпространств Крылова равна $n$. Соответственно, $n$ ортогональных векторов $h_0, h_1, \ldots, h_{n-1}$ образуют в $K_n$ базис. Поскольку градиент в любой порождаемой методом первого порядка точке принадлежит соответствующему подпространству Крылова, $\nabla J(q_n)$ раскладывается по этому базису. По следствию из леммы 1, это обозначает, что $q_n=q^*$. Но в таком случае утверждение теоремы тоже выполняется.

QED. $ $

Важность доказанной теоремы следует из того, что все моментные методы, включая различные варианты метода сопряжённых градиентов, лежат в классе методов, работающих в подпространствах Крылова. Теорема показывает, что $\infty$-моментный метод минимальных ошибок обеспечивает наилучшую сходимость по аргументу в этом классе.

{\bf Теорема 2. О сходимости по аргументу метода градиентного спуска.}

    Если некорректно поставленная задача (\ref{Aq=f}) в гильбертовом пространстве $H$ имеет единственное решение $q^*$ и оператор $A_0$ компактен, то простой градиентный спуск $$q_{k+1}=q_k-\alpha\nabla J(q_k)$$ (при $\alpha\leq\frac{1}{L}$) сходится по аргументу к точному решению, то есть $$\lim\limits_{k\rightarrow\infty}||q_k-q^*||^2 = 0.$$

{\bf Доказательство.}
    Если оператор $A_0$ компактен, то сопряжённый ему $A^*$ тоже компактен (\cite{kolmogorov2004}, гл.4, $\S6$, теорема 3). Тогда оператор $B=A^*A_0$ --- самосопряжённый компактный, и по теореме Гильберта---Шмидта (\cite{kolmogorov2004}, гл.4, $\S6$, теорема 5) существуют монотонно стремящаяся к нулю последовательность его собственных значений $\{\lambda_n\}_{n=1}^{\infty}$ и последовательность соответствующих им собственных векторов $\{w_n\}_{n=1}^{\infty}$ такие, что произвольный вектор $x\in H$ может быть представлен в виде $x=\sum\limits_{j=1}^{\infty}\alpha_j w_j + x'$, где $x'\in\mathrm{Ker}B$, и $Bx=\sum\limits_{j=1}^{\infty}\alpha_j\lambda_j w_j$.

    Поскольку уравнение $Aq=f$ имеет единственное решение, оператор $A_0$ (а значит и $A^*$) инъективен, и ядро тривиально. Следовательно $x'=0$, и $\{w_n\}_{n=1}^{\infty}$ --- базис в $H$.

    Пусть $q\in H$ --- произвольный вектор из $H$. Разложим $q-q^*$ по базису: $$q-q^*=\sum\limits_{n=1}^{\infty}\xi_nw_n.$$
    Тогда можно вычислить градиент $$\nabla J(q) = B(q-q^*)=\sum\limits_{n=1}^{\infty}\lambda_n\xi_nw_n$$ и функционал $$J(q)=\frac{1}{2}||A_0(q-q^*)||^2 = \frac{1}{2}\langle q-q^*, B(q-q^*)\rangle = \frac{1}{2}\sum\limits_{n=1}^{\infty}\lambda_n\xi_n^2.$$

    Применим в точке $q$ градиентный спуск с шагом $\alpha\leq\frac{1}{L}=\frac{1}{\lambda_1}$:
    $$q_{new} = q - \alpha\nabla J(q) = q^* + \sum\limits_{n=1}^{\infty}\xi_n\left(1-\alpha\lambda_n\right)w_n.$$
    То есть, каждая компонента просто умножается на своё число из полуинтервала $[0;1)$, потому что $1>1-\alpha\lambda_n\geq1-\frac{\lambda_n}{\lambda_1}\geq 0$.

    Пусть $q_0=q^*+\sum\limits_{n=1}^{\infty}\xi_nw_n$ --- это начальное приближение. Тогда после $k$ шагов получается приближение $$q_k=q^* + \sum\limits_{n=1}^{\infty}\xi_n\left(1-\alpha\lambda_n\right)^kw_n.$$
    Расстояние от него до точного решения вычисляется стандартно:
    $$||q_k-q^*||^2 = \sum\limits_{n=1}^{\infty}\xi_n^2\left(1-\alpha\lambda_n\right)^{2k}.$$
    Этот ряд сходится, в том числе и при $k=0$, поскольку $\xi_k$ --- коэффициенты Фурье некоторого элемента $H$.

    Из сходимости ряда следует, что $\forall\varepsilon>0\ \exists N=N(\varepsilon)\in\mathbb{N}: \sum\limits_{n=N+1}^{\infty}\xi_n^2 <\frac{\varepsilon}{2}.$

    Для всех $n\in\{1,\ldots,N\}$ существуют такие $K_n\in\mathbb{N}$, что при $k\geq K_n$ $\xi_n^2\left(1-\alpha\lambda_n\right)^{2k} < \frac{\varepsilon}{2N}$. Обозначим $K=\max\{K_1,\ldots,K_N\}$. По определению, при фиксированной последовательности $\xi_n$ (коэффициентов разложения начального приближения) оно зависит только от выбора $\varepsilon$.

    Соответственно, $\forall\varepsilon>0\ \exists K=K(\varepsilon)\in\mathbb{N}: \forall k\geq K\ \sum\limits_{n=1}^{\infty}\xi_n^2\left(1-\alpha\lambda_n\right)^{2k} = \\$$ $$=\sum\limits_{n=1}^{N}\xi_n^2\left(1-\alpha\lambda_n\right)^{2k} + \sum\limits_{n=N+1}^{\infty}\xi_n^2\left(1-\alpha\lambda_n\right)^{2k} < N\cdot\frac{\varepsilon}{2N} + \frac{\varepsilon}{2} = \varepsilon$.

    А это и есть определение нулевого предела.

QED. $ $

{\bf Теорема 3. О сходимости по аргументу $\infty$-моментного метода минимальных ошибок.}

    В условиях теоремы 2, $\infty$-моментный метод минимальных ошибок сходится по аргументу к точному решению, то есть $$\lim\limits_{k\rightarrow\infty}||q_k-q^*||^2 = 0.$$
    
{\bf Доказательство.}
    Пусть $q_0$ --- начальное приближение; $\{q_k\}_{k=0}^{\infty}$ --- последовательность точек, порождаемых $\infty$-моментным методом минимальных ошибок; $\{\tilde{q}_k\}_{k=0}^{\infty}$ --- последовательность точек, порождаемых простейшим градиентным спуском. По теореме 1, $||q_k-q^*||\leq ||\tilde{q}_k-q^*||$ при всех натуральных $k$. По теореме 2, $||\tilde{q}_k-q^*||\rightarrow 0$ при $k\rightarrow\infty$. Из теоремы о двух милиционерах получаем требуемое.

QED. $ $

Строгое доказательство этого простого следствия необходимо. Ведь в конечномерном пространстве максимальная по включению линейно независимая система образует базис (что гарантирует достижение решения рассматриваемым методом), а в бесконечномерном пространстве гипотетически возможна ситуация, когда $q_0-q^*$ не раскладывается по системе $\{g,Bg,\ldots\}$, и достижение $q^*$ оказывается невозможным. Данная теорема гарантирует сходимость к решению, если оно единственно. 

Если же решений больше одного (для квадратичной задачи это обозначает существование линейного многообразия решений), то метод будет сходиться к одному из них. 

{\bf Теорема 4. О скорости сходимости по функционалу $m$-моментного метода минимальных ошибок.}

$m$-моментный метод минимальных ошибок ($m\geq 0$) сходится по функционалу со сверхлинейной скоростью, то есть $$J(q_k)=o\left(\frac{1}{k}\right), k\rightarrow\infty.$$

{\bf Доказательство.}
    Сначала проведём вычисления для метода минимальных ошибок ($m=0$). 
$$J(q_k) = \frac{1}{2}\alpha_k||s_k||^2 = \frac{1}{2\alpha_k}||q_{k+1}-q_k||^2 = \frac{||\nabla J(q_k)||^2}{4J(q_k)}||q_{k+1}-q_k||^2\Rightarrow$$
$$\Rightarrow J(q_k) = \frac{1}{2}||\nabla J(q_k)||\cdot||q_{k+1}-q_k||.$$

Поскольку для функционала с липшицевым градиентом выполняется неравенство $$J(q)\geq \frac{||\nabla J(q)||^2}{2L},$$ можно получить следующую оценку:
$$J(q_k)\leq \frac{1}{2}\sqrt{2L J(q_k)}\cdot||q_{k+1}-q_k||\Rightarrow J(q_k)\leq \frac{L}{2}||q_{k+1}-q_k||^2,$$
откуда и вытекает требуемая оценка скорости сходимости: из сходимости ряда (\ref{limNorm2}) следует стремление общего члена к нулю быстрее $\frac{1}{k}$.

При $m\geq 1$ вычисления очень похожи:
$$J(q_k) = \frac{1}{2}\alpha_k||s_k||^2 = \frac{1}{2\alpha_k}||q_{k+1}-q_k||^2 = \frac{||\nabla J(q_k)||^2\sin^2{\varphi_k}}{4J(q_k)}||q_{k+1}-q_k||^2\Rightarrow$$
$$\Rightarrow J(q_k) = \frac{1}{2}||\nabla J(q_k)||\cdot||q_{k+1}-q_k||\sin{\varphi_k}.$$
Подставляя оценку нормы градиента, получаем:
$$J(q_k)\leq \frac{L}{2}||q_{k+1}-q_k||^2\sin^2{\varphi_k}.$$

Это обозначает сходимость не хуже метода минимальных ошибок.

QED. $ $

{\bf Теорема 5. О скорости сходимости по аргументу $m$-моментного метода минимальных ошибок.}

Метод минимальных ошибок для сильно выпуклого квадратичного функционала сходится по аргументу со скоростью геометрической прогрессии, то есть $$||q_{k+1}-q^*||^2\leq \left(1-\frac{\mu}{L}\right)\cdot||q_k-q^*||^2.$$
Для $m$-моментного метода минимальных ошибок оценка сходимости улучшается: $$||q_{k+1}-q^*||^2\leq \left(1-\frac{\mu}{L\sin^2{\varphi_k}}\right)\cdot||q_k-q^*||^2.$$

{\bf Доказательство.}
Воспользуемся оценками, полученными при доказательстве прошлой теоремы, и свойством сильно выпуклого функционала $$J(q)\geq \frac{\mu||q-q^*||^2}{2}.$$

Для простого метода минимальных ошибок:
$$||q_{k+1}-q_k||^2\geq \frac{2J(q_k)}{L}\geq \frac{\mu||q_k-q^*||^2}{L}\Rightarrow$$ $$\Rightarrow ||q_{k+1}-q^*||^2 = ||q_k-q^*||^2 - ||q_{k+1}-q_k||^2\leq \left(1-\frac{\mu}{L}\right)\cdot||q_k-q^*||^2,$$
что и требовалось доказать.

Эта оценка соответствует обычным оценкам сходимости для градиентных методов.

Для $m$-моментного метода минимальных ошибок:
$$||q_{k+1}-q_k||^2\geq \frac{2J(q_k)}{L\sin^2{\varphi_k}}\geq \frac{\mu||q_k-q^*||^2}{L\sin^2{\varphi_k}}\Rightarrow$$ $$\Rightarrow ||q_{k+1}-q^*||^2 = ||q_k-q^*||^2 - ||q_{k+1}-q_k||^2\leq \left(1-\frac{\mu}{L\sin^2{\varphi_k}}\right)\cdot||q_k-q^*||^2,$$
что и требовалось доказать.

Эта оценка лучше предыдущей. Убывание расстояния до точного решения зависит от текущей величины угла между антиградиентом и его проекцией на линейную оболочку векторов $m$ предыдущих шагов. Если этот угол близок к $\frac{\pi}{2}$, то для данного шага оценка соответствует обычной оценке для градиентных методов. Если же угол мал, то наблюдается очень быстрое убывание.

QED. $ $

Получены оценки сходимости в сильно выпуклом случае; в конечномерном пространстве доказана сходимость за конечное число шагов. Возникает вопрос, нет ли хороших глобальных оценок скорости сходимости по аргументу и для не сильно выпуклых задач оптимизации, которые возникают при решении некорректных задач? Отрицательный ответ на этот вопрос даёт следующая теорема.

{\bf Теорема 6. О неравномерной сходимости по аргументу; или о существовании начального приближения, при котором сходимость по аргументу оказывается сколь угодно медленной.}

    В условиях теоремы 2, для любых наперёд заданных $\varepsilon\in (0, 1)$ и $N\in\mathbb{N}$ существует такое начальное приближение $q_0$, что $||q_0-q^*||=1$, а для любого метода, работающего в подпространствах Крылова, в том числе и для $\infty$-моментного метода минимальных ошибок, $||q_N-q^*||^2>\varepsilon$.
    
{\bf Доказательство.}
    Зафиксируем произвольное натуральное число $N$ и $\varepsilon\in (0, 1)$.  
    
    $\nabla J(q_0)=B(q_0-q^*)$, поэтому $K_N=Lin\left(\{B^k(q_0-q^*)\}_{k=1}^N\right)$. При доказательстве теоремы 2 получен ортонормированный базис из собственных векторов $B$, будем использовать разложение по нему:
    $$q_0-q^*=\sum\limits_{n=1}^{\infty}\xi_nw_n;\quad B^k(q_0-q^*)=\sum\limits_{n=1}^{\infty}\lambda_n^k\xi_nw_n.$$

    Будем считать, что каждому $\lambda_n$ соответствует только один собственный вектор, то есть последовательность собственных значений --- строго убывающая. Если это не так, и каким-то собственным значениям соответствует несколько собственных векторов, то будем рассматривать лишь такие начальные приближения $q_0$, для которых не более одного коэффициента $\xi_n$ при векторе из такого набора отличается от нуля. Или просто (аналогично методу вспомогательного угла из школьной тригонометрии) представим линейную комбинацию собственных векторов, соответствующих одному и тому же собственному значению, в виде одного такого вектора, умноженного на число: $\sum\limits_{i=1}^k\xi_{n_i}w_{n_i}=r\sum\limits_{i=1}^k\frac{\xi_{n_i}}{r}w_{n_i}$, где $r=\sqrt{\sum\limits_{i=1}^k\xi_{n_i}^2}$.

    Любой метод, работающий в подпространствах Крылова, порождает на шаге $N$ точку $q_N=q_0+\sum\limits_{i=1}^N\eta_iB^i (q_0-q^*)$ для некоторого набора чисел $\eta_1,\ldots,\eta_N$. Преобразуем расстояние до точного решения: $$||q_N-q^*||^2=||q_0-q^*+\sum\limits_{i=1}^N\eta_iB^i (q_0-q^*)||^2 = \left|\left|\sum_{n=1}^{\infty} \xi_n\left(1+\sum\limits_{i=1}^{N}\eta_i\lambda_n^i\right)w_n\right|\right|^2=$$ $$=\sum\limits_{n=1}^{\infty}\xi_n^2\left(1+\sum\limits_{i=1}^{N}\eta_i\lambda_n^i\right)^2=:\Psi(\eta; q_0).$$

    $\Psi(\eta; q_0)$ --- квадратичная функция от $N$-мерного вектора $\eta=\left(\eta_1,\ldots,\eta_N\right)^T$, зависящая от параметра $q_0$. Наша цель --- подобрать такие $\xi_n$, что $\sum\limits_{n=1}^{\infty} \xi_n^2=1$, а $\min\limits_{\eta\in\mathbb{R}^n}\Psi\left(\eta; \sum\limits_{n=1}^{\infty}\xi_n w_n\right)>\varepsilon$. Тогда утверждение теоремы будет выполнено.

    Заметим: если количество ненулевых $\xi_n$ не превосходит $N$, то $\min\limits_{\eta\in\mathbb{R}^n}\Psi\left(\eta; \sum\limits_{n=1}^{\infty}\xi_n w_n\right)=0$. Действительно, для обнуления $\Psi\left(\eta; \sum\limits_{n=1}^{\infty}\xi_n w_n\right)$ необходимо и достаточно выполнение равенства $1+\sum\limits_{i=1}^{N}\eta_i\lambda_n^i=0$ для всех таких $n$, что $\xi_n\not=0$ (обычное условие стационарности через частные производные, естественно, при этом выполняется). Если количество уравнений в этой системе оказалось меньше $N$, дополним её до этого количества и другими номерами. Определитель получившейся системы --- это известный определитель Вандермонда (помноженный на произведение $\lambda_n$), а он не равен нулю при различных положительных $\lambda_n$. Поэтому решение существует. 

    Возьмём $\xi=\xi_1=\xi_2=\ldots=\xi_N=\sqrt{\frac{1-\varepsilon}{2N}}$, $\xi_M=\sqrt{\frac{1+\varepsilon}{2}}$, где $M>N$ --- число, которое будет выбрано позднее; все остальные $\xi_n$ равны нулю. Условие нормировки выполняется. Обозначим $\tilde{\eta}$ --- набор, при котором $\Psi\left(\tilde{\eta}; \sum\limits_{n=1}^N\xi w_n\right)=0$. По доказанному в предыдущем абзаце, он существует.

    Запишем необходимое условие минимальности $\Psi\left(\eta; \sum\limits_{n=1}^N\xi w_n + \xi_Mw_M\right)$ в частных производных:
    $$\frac{1}{2}\frac{\partial}{\partial\eta_j}\Psi\left(\eta; \sum\limits_{n=1}^N\xi w_n + \xi_Mw_M\right) = \frac{1}{2}\frac{\partial}{\partial\eta_j}\left(\xi^2\sum\limits_{n=1}^N\left(1+\sum\limits_{i=1}^{N}\eta_i\lambda_n^i\right)^2 + \xi_M^2\left(1+\sum\limits_{i=1}^{N}\eta_i\lambda_M^i\right)^2\right) =$$
    
    \begin{equation}
    \label{Psi_minimum}
        = \xi^2\sum\limits_{n=1}^N\left(1+\sum\limits_{i=1}^{N}\eta_i\lambda_n^i\right)\lambda_n^j + \xi_M^2\left(1+\sum\limits_{i=1}^{N}\eta_i\lambda_M^i\right)\lambda_M^j = 0 \text{ при } 1\leq j\leq N.
    \end{equation}

    Обозначим вектор $\left(1, \lambda_n,\ldots,\lambda_n^{N-1}\right)^T=:x_n$. Тогда система уравнений (\ref{Psi_minimum}) принимает вид 
    $$\left(\sum\limits_{n=1}^N\xi^2\lambda_n^2x_nx_n^T + \xi_M^2\lambda_M^2x_Mx_M^T\right)\eta = -\xi^2\sum\limits_{n=1}^N\lambda_nx_n - \xi_M^2\lambda_Mx_M
    $$

    Подставим значения $\xi$ и $\xi_M$, а затем домножим левую и правую части на одно и то же число:
    
    \begin{equation}
    \label{Psi_minimum_norm}\left(\sum\limits_{n=1}^N\lambda_n^2x_nx_n^T + N\frac{1+\varepsilon}{1-\varepsilon}\lambda_M^2x_Mx_M^T\right)\eta = -\sum\limits_{n=1}^N\lambda_nx_n - N\frac{1+\varepsilon}{1-\varepsilon}\lambda_Mx_M.\end{equation}

    Матрица $L=\sum\limits_{n=1}^N\lambda_n^2x_nx_n^T$ --- положительно определённая. Действительно: предположим, что для некоторого $y\not=\overline{0}$ $y^TLy=0$. Тогда $\sum\limits_{n=1}^N\lambda_n^2\langle x_n, y\rangle^2 = 0$. Поскольку все $\lambda_n>0$, это возможно, только если $\langle x_n, y\rangle=0$ при всех $n$. Но векторы $x_n$, как уже было замечено, --- это $N$ линейно независимых векторов в $\mathbb{R}^N$. Поэтому ненулевой $y\in\mathbb{R}^N$ не может быть ортогонален им всем.

    Значит, система (\ref{Psi_minimum_norm}) имеет единственное решение, поскольку её матрица --- это сумма положительно определённой $L$ и неотрицательно определённой матрицы, заданной последним слагаемым. Обозначим это решение $\hat{\eta}(M)$.

    Поскольку набор $\tilde{\eta}$ обнуляет все остальные слагаемые в определении, $$\Psi\left(\tilde{\eta}; \sum\limits_{n=1}^N\xi w_n+\xi_Mw_M\right)=\left|\left|\xi_M\left(1+\sum\limits_{i=1}^n\tilde{\eta}_i\lambda_M^i\right)w_i\right|\right|^2=\frac{1+\varepsilon}{2}\left(1+\sum\limits_{i=1}^n\tilde{\eta}_i\lambda_M^i\right)^2.$$

    При $M\rightarrow\infty$ $\lambda_M\rightarrow 0$, поэтому выражение в скобке стремится к $1$. Соответственно, при достаточно больших $M$ выполняется оценка 
    
    \begin{equation}
    \label{Psi1}
        \Psi\left(\tilde{\eta}; \sum\limits_{n=1}^N\xi w_n+\xi_Mw_M\right) > \dfrac{1+3\varepsilon}{4}
    \end{equation}

    Известно, что при положительно определённой матрице решение системы линейных алгебраических условий непрерывно зависит от левой и правой частей: если $\delta L\rightarrow 0$ и $\delta b\rightarrow 0$ по некоторой норме (например, по евклидовой; хотя в конечномерном пространстве, в котором решается система (\ref{Psi_minimum_norm}), они все эквивалентны), то решение системы $(L+\delta L)\eta=b+\delta b$ стремится к решению $L\eta=b$ по той же норме. 

    Обозначим слагаемые формулы (\ref{Psi_minimum_norm}), содержащие номер $M$, за $\delta L$ и $\delta b$ соответственно. Заметим также, что при $M\rightarrow\infty$ $\lambda_M\rightarrow 0$, $||x_M||\rightarrow 1$. Значит, $||\hat{\eta}(M)-\tilde{\eta}||\rightarrow 0$, и, поскольку $\Psi(\eta;q)$ --- непрерывная функция, $\left|\Psi\left(\hat{\eta}(M); \sum\limits_{n=1}^N\xi w_n\right) - \Psi\left(\tilde{\eta}; \sum\limits_{n=1}^N\xi w_n\right)\right| \rightarrow 0$.

    $\left|\Psi\left(\eta; \sum\limits_{n=1}^N\xi w_n + \xi_M w_M\right) - \Psi\left(\eta; \sum\limits_{n=1}^N\xi w_n\right)\right| = \frac{1+\varepsilon}{2}\left(1+\sum\limits_{i=1}^n\eta_i\lambda_M^i\right)^2$. При $M\rightarrow\infty$ эта величина стремится к $\frac{1+\varepsilon}{2}$, даже если $\eta$ не фиксирован, а зависит от номера $M$ и лишь ограничен по норме.

    Значит, $\Psi\left(\hat{\eta}(M); \sum\limits_{n=1}^N\xi w_n+\xi_Mw_M\right) - \Psi\left(\tilde{\eta}; \sum\limits_{n=1}^N\xi w_n+\xi_Mw_M\right) = \Delta_1 + \Delta_2 - \Delta_3$, где
    $$\Delta_1=\Psi\left(\hat{\eta}(M); \sum\limits_{n=1}^N\xi w_n\right) - \Psi\left(\tilde{\eta}; \sum\limits_{n=1}^N\xi w_n\right)\rightarrow 0\text{ при } M\rightarrow\infty;$$ $$\Delta_2= \Psi\left(\hat{\eta}(M); \sum\limits_{n=1}^N\xi w_n+\xi_Mw_M\right) - \Psi\left(\hat{\eta}(M); \sum\limits_{n=1}^N\xi w_n\right)\rightarrow \dfrac{1+\varepsilon}{2}\text{ при } M\rightarrow\infty;$$ $$\Delta_3=\Psi\left(\tilde{\eta}; \sum\limits_{n=1}^N\xi w_n+\xi_Mw_M\right) - \Psi\left(\tilde{\eta}; \sum\limits_{n=1}^N\xi w_n\right)\rightarrow \dfrac{1+\varepsilon}{2}\text{ при } M\rightarrow\infty.$$

    Соответственно, $\Psi\left(\hat{\eta}(M); \sum\limits_{n=1}^N\xi w_n+\xi_Mw_M\right) - \Psi\left(\tilde{\eta}; \sum\limits_{n=1}^N\xi w_n+\xi_Mw_M\right)\rightarrow 0$ при $M\rightarrow\infty$. Поэтому, начиная с некоторого $M$, выполняется оценка 
    \begin{equation}
    \label{Psi2}
        \left|\Psi\left(\hat{\eta}(M); \sum\limits_{n=1}^N\xi w_n+\xi_Mw_M\right) - \Psi\left(\tilde{\eta}; \sum\limits_{n=1}^N\xi w_n+\xi_Mw_M\right)\right|<\dfrac{1-\varepsilon}{4}
    \end{equation}

    Выберем такое большое $M$, чтобы выполнялись обе оценки (\ref{Psi1}) и (\ref{Psi2}). Тогда для $\hat{\eta}(M)$ --- точки минимума функции $\Psi$ --- выполняется следующая оценка: \\$\Psi\left(\hat{\eta}(M); \sum\limits_{n=1}^N\xi w_n+\xi_Mw_M\right) > \dfrac{1+3\varepsilon}{4} - \dfrac{1-\varepsilon}{4}=\varepsilon$, что и требовалось доказать.

QED. $ $

Эта теорема показывает, что, хотя для любого начального приближения и достигается сходимость $\infty$-моментного метода минимальных ошибок (то есть, стремление к нулю расстояния до точного решения), равномерной для всех точек оценки этой сходимости не существует. За фиксированное количество шагов сокращение расстояния может быть сколь угодно слабым.

Важно отметить, что доказанное свойство сколь угодно медленной сходимости относится не к $\infty$-моментному методу минимальных ошибок: на нём используемый в доказательстве минимум реализуется, а для любого другого метода первого порядка сходимость будет ещё медленнее (теорема 1). Оно относится именно к бесконечномерному пространству, в котором решаются задачи оптимизации.

\vspace{12pt}
\centerline{5. Экспериментальная проверка построенного метода}
\vspace{12pt}
\centerline{\it 5.1. Решение задачи Коши для уравнения Гельмгольца}
\vspace{12pt}

Уравнение Гельмгольца относится к эллиптическому типу. Для таких уравнений корректными являются краевые задачи, а некорректными --- начально-краевые. В качестве иллюстрации для сравнения методов использована модельная задача \begin{equation}
\begin{cases}u_{xx}+u_{yy} + \kappa^2 u = r(x,y), & (x,y)\in \Omega = (0,1)\times (0,1)\\u|_{x=0}=f(y), & y\in[0,1]\\u_x|_{x=0}=g(y), & y\in[0,1]\\u|_{y=0}=u|_{y=1}=0, & x\in[0,1]\end{cases}
\label{helmholtz}
\end{equation} с параметрами $r(x, y)=-x(2-y+y^2)$, $f(y)=0$, $g(y)=y-y^2$. 

Начально-краевая задача для уравнения Гельмгольца некорректна, но корректна краевая задача:
\begin{equation}
\begin{cases}\Delta u + \kappa^2 u = r(x, y), & (x,y)\in \Omega = (0,1)\times (0,1)\\u_x|_{x=0} = g(y), & y\in[0,1]\\u|_{y=0} = u|_{y=1} = 0, & x\in[0,1]\\u|_{x=1}=q(y), & y\in[0,1]\\\end{cases}.
\label{helmholtz_correct}
\end{equation}
Эта задача отличается от исходной лишь заменой условия $u|_{x=0}=f(y)$ на $u|_{x=1}=q(y)$. Граница $x=0$ --- наблюдаемая, $x=1$ --- ненаблюдаемая. Требуется найти $q(y)$. Для данной задачи известно точное решение: $q(y)=y-y^2$.

На пространстве непрерывно дифференцируемых функций со стандартным скалярным произведением $$H = \{q\in C^1[0,1]\cap L^2[0,1]: q(0)=q(1)=0\}$$ определим оператор $A: H\rightarrow H$ следующим образом:

\begin{equation}
(Aq)(y)=u(0,y)\quad \forall y\in(0,1),
\label{helmholtzA}
\end{equation}
где $u$ --- решение корректной краевой задачи (\ref{helmholtz_correct}).

Эта постановка в виде операторного уравнения приведена в статье \cite{pletnev2022helmholtz}; там же вычислен градиент функционала.

\begin{table}[H]
\begin{center}
\begin{tabular}{|c|c|c|}\hline
    $ $ & $||q_n-q^*||$ & $J(q_n)$ \\
   \hline
   ММО, $m=1$ & $6.14\times 10^{-4}$ & $7.38\times 10^{-24}$ \\
   \hline
   ММО, $m=2$ & $6.14\times 10^{-4}$ & $3.32\times 10^{-25}$ \\
   \hline
   ММО, $m=5$ & $6.14\times 10^{-4}$ & $1.53\times 10^{-25}$ \\
   \hline
   ММО, $m=\infty$ & $6.14\times 10^{-4}$ & $2.19\times 10^{-25}$ \\
   \hline
   Сопряжённые градиенты & $1.58\times 10^{-3}$ & $1.29\times 10^{-19}$ \\
   \hline
   Тяжёлый шарик & $1.58\times 10^{-3}$ & $1.36\times 10^{-19}$ \\
   \hline
   Шаг Поляка & $6.94\times 10^{-3}$ & $8.30\times 10^{-13}$ \\
   \hline
   STM & $1.56\times 10^{-3}$ & $1.24\times 10^{-19}$ \\
   \hline
\end{tabular}
\end{center}
\caption{Результаты работы методов при решении задачи (\ref{helmholtz}) с использованием сетки с шагом $h=0.01$. Начальное значение функционала $J(q_0) = 1.77\cdot 10^{-4}$, начальное расстояние до точного решения $||q_0-q^*||=0.183$.}
\label{tab:table4.1}
\end{table}

\begin{figure}[H]
	\begin{minipage}[h]{0.5\linewidth}
		\center{\includegraphics[width=0.9\linewidth]{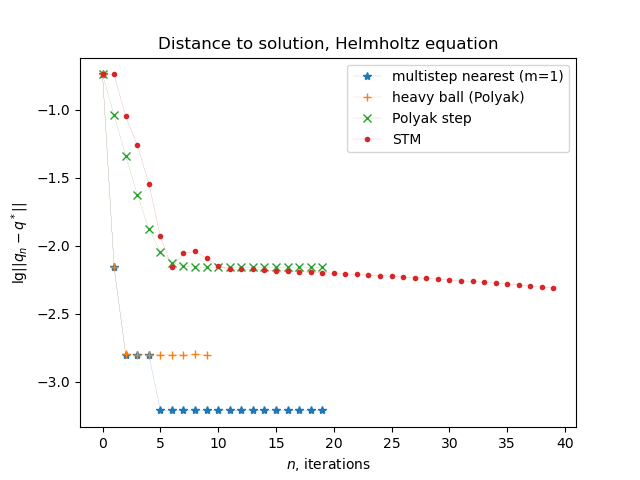}}
	\end{minipage}
	\hfill
	\begin{minipage}[h]{0.5\linewidth}
		\center{\includegraphics[width=0.9\linewidth]{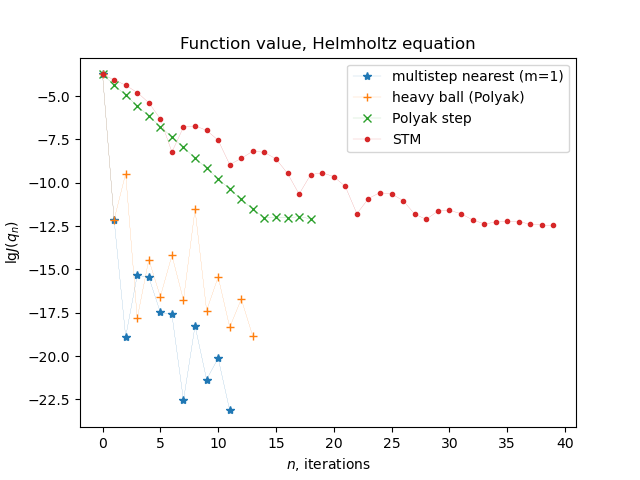}}
	\end{minipage}
	\caption{Сравнение $1$-моментного ММО, адаптивного тяжёлого шарика, градиентного спуска с шагом Поляка и метода подобных треугольников в применении к решению задачи (\ref{helmholtz}) с использованием сетки с шагом $h=0.01$.}
	\label{ris:image4.1}
\end{figure}

Рисунок \ref{ris:image4.1} показывает, что $1$-моментный метод минимальных ошибок позволяет достичь меньшей невязки, чем адаптивный метод тяжёлого шарика из статьи \cite{goujaud2022quadratic}, градиентный спуск с шагом Поляка и метод подобных треугольников, и по аргументу, и по функционалу. При этом сходимость по функционалу не является монотонной. 

\begin{figure}[H]
	\begin{minipage}[h]{0.5\linewidth}
		\center{\includegraphics[width=0.9\linewidth]{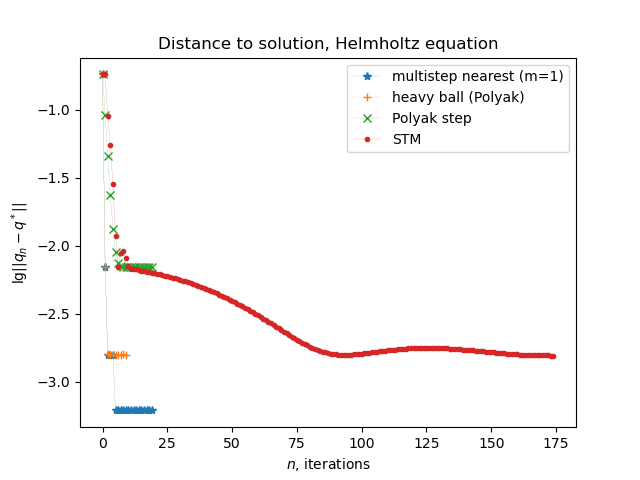}}
	\end{minipage}
	\hfill
	\begin{minipage}[h]{0.5\linewidth}
		\center{\includegraphics[width=0.9\linewidth]{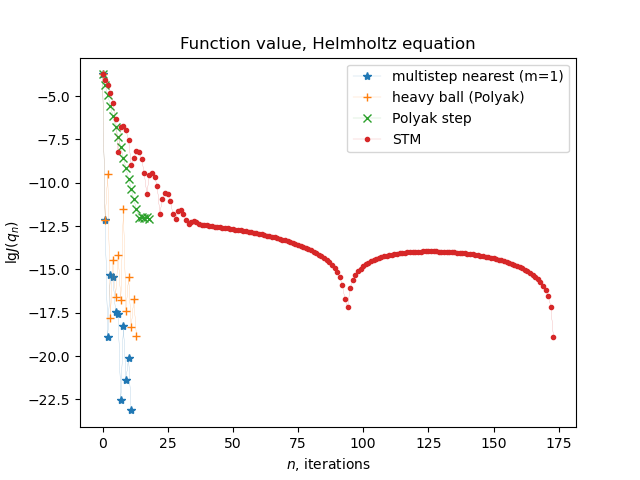}}
	\end{minipage}
	\caption{Сравнение $1$-моментного ММО, адаптивного тяжёлого шарика, градиентного спуска с шагом Поляка и метода подобных треугольников в применении к решению задачи (\ref{helmholtz}) с использованием сетки с шагом $h=0.01$. Отмечено $175$ точек для метода подобных треугольников.}
	\label{ris:image4.1.1}
\end{figure}

Увеличение количества итераций не приводит к заметным изменениям результатов, кроме метода подобных треугольников. Более длинный отрезок приведён на рисунке \ref{ris:image4.1.1}

\begin{figure}[H]
	\begin{minipage}[h]{0.5\linewidth}
		\center{\includegraphics[width=0.9\linewidth]{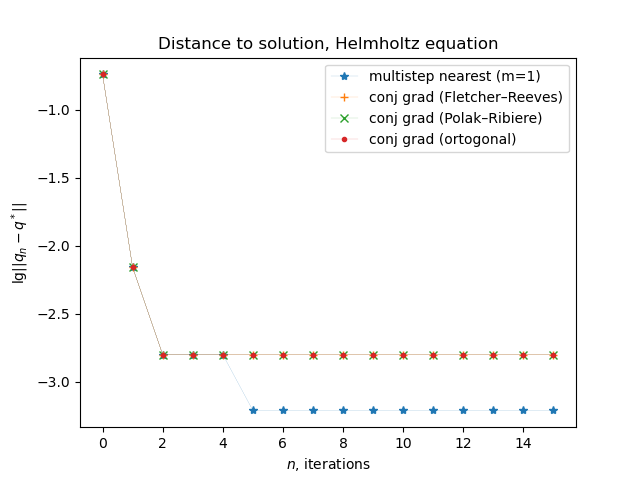}}
	\end{minipage}
	\hfill
	\begin{minipage}[h]{0.5\linewidth}
		\center{\includegraphics[width=0.9\linewidth]{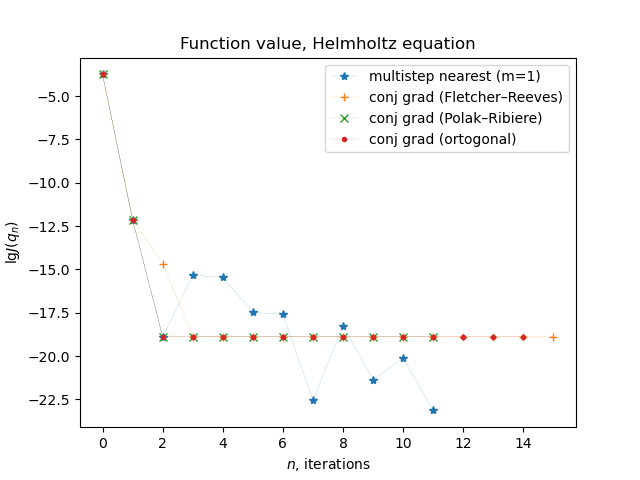}}
	\end{minipage}
	\caption{Сравнение $1$-моментного ММО и различных вариантов метода сопряжённых градиентов в применении к решению задачи (\ref{helmholtz}) с использованием сетки с шагом $h=0.01$.}
	\label{ris:image4.2}
\end{figure}

Рисунок \ref{ris:image4.2} показывает, что различные варианты метода сопряжённых градиентов дают на рассматриваемой задаче примерно одинаковые результаты, а $1$-моментный метод минимальных ошибок их превосходит. Интересно, что достигаемое значение функционала для него тоже оказывается существенно меньше.

\begin{figure}[H]
	\begin{minipage}[h]{0.5\linewidth}
		\center{\includegraphics[width=0.9\linewidth]{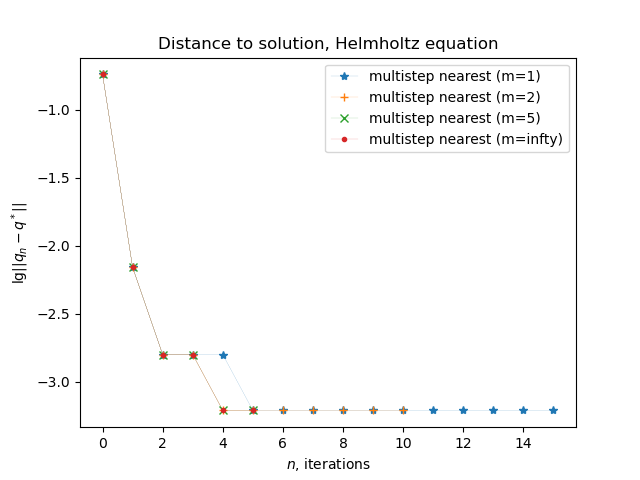}}
	\end{minipage}
	\hfill
	\begin{minipage}[h]{0.5\linewidth}
		\center{\includegraphics[width=0.9\linewidth]{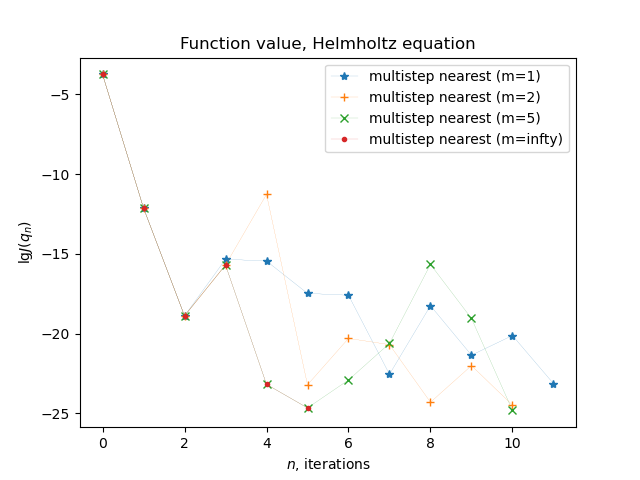}}
	\end{minipage}
	\caption{Сравнение $m$-моментных ММО с $m\in\{1, 2, 5, \infty\}$ в применении к решению задачи (\ref{helmholtz}) с использованием сетки с шагом $h=0.01$.}
	\label{ris:image4.3}
\end{figure}

На рисунке \ref{ris:image4.3} сравниваются $m$-моментные методы минимальных ошибок при разных $m$. Теоретически, чем больше $m$, тем лучше должна быть сходимость по аргументу. Однако на данной задаче такой эффект не наблюдается: все методы дают очень близкие результаты.

Эксперименты показывают, что градиентные методы оптимизации позволяют эффективно решать начально-краевую задачу для уравнения Гельмгольца.

Это согласуется с теоретическими результатами: оператор задачи компактен, поэтому по теореме 2 сходится метод простейшего градиентного спуска; соответственно, рассмотренные методы тоже должны сходиться. В частности, сходимость $m$-моментного метода минимальных ошибок установлена теоремой 3.

Теорема 1 об оптимальности $\infty$-моментного метода минимальных ошибок нашла своё подтверждение: новые методы сходятся значительно лучше ранее существовавших. Ускоренные методы достигли лучших результатов, чем неускоренный (градиентный спуск с шагом Поляка).

Как показано в статье \cite{pletnev2022helmholtz}, собственные значения связанного с задачей самосопряжённого оператора $A^*A_0$ имеют вид $\lambda_n=\frac{1}{\mathrm{ch}^2 \sqrt{\pi^2 n^2-\kappa^2}}$. 
Поскольку при $t\rightarrow\infty$ $\mathrm{ch} t\sim \frac{e^t}{2}$, а $\sqrt{\pi^2n^2-\kappa^2}=\pi n - \frac{\kappa^2}{2\pi n} + o\left(\frac{1}{n}\right)$, $n\rightarrow\infty$, собственные значения стремятся к нулю эквивалентно геометрической прогрессии со знаменателем $e^{2\pi}$.

Тот факт, что последовательность собственных значений связанного с задачей самосопряжённого оператора стремится к нулю не слишком быстро, может быть причиной быстрого достижения градиентными методами оптимизации точек, близких к точному решению. Это связано с тем, что такая задача хорошо приближается конечномерной, если в рядах Фурье для $f(y)$ и $g(y)$ отбросить слагаемые, начиная с некоторого номера. С одной стороны, остаток ряда получается достаточно маленьким; с другой стороны, константа сильной выпуклости (минимальное собственное значение, соответствующее оставшимся слагаемым) оказывается не слишком близкой к нулю.

\vspace{12pt}
\centerline{\it 5.2. Решение ретроспективной задачи Коши для уравнения теплопроводности}
\vspace{12pt}

В статье \cite{pletnev2023component} показано, что ретроспективная задача Коши для уравнения теплопроводности является некорректной. Это вполне соответствует физическим соображениям. В качестве иллюстрации для сравнения методов использована модельная задача 
\begin{equation}
\begin{cases}u_t - \kappa^2 \Delta_x u = 0, & (x,t)\in \Omega = \Pi\times (0,1)\\u|_{x\in\partial \Pi}=0, & t\in[0,1]\\u|_{t=1} = f(x), & x\in\Pi\end{cases}
\label{heat3Du}
\end{equation}
($\Pi=(0,1)^3$) с переменным коэффициентом теплопроводности

$$\kappa(x) = 
\begin{cases}
    \kappa_{max}, \text{ если } 0.4<x_1, x_2, x_3<0.6; &\\
    \frac{\kappa_{max}}{5}, \text{ иначе.}
\end{cases},$$
где $\kappa_{max}\in\{0.4,0.6\}$, точным решением $q(x)=\sin{2\pi x_1}\cdot\sin^2{2\pi x_2}\cdot\sin^3{2\pi x_3}$ и условием на наблюдаемой границе $f(x)=(Aq)(x)$.

Корректная задача с условием в начальный момент времени
\begin{equation}
\begin{cases}u_t - \kappa^2 \Delta_x u = 0, & (x,t)\in \Omega = \Pi\times (0,1)\\u|_{x\in\partial \Pi}=0, & t\in[0,1]\\u|_{t=0} = q(x), & x\in\Pi\end{cases}
\label{heat3Du_correct}
\end{equation} отличается от исходной лишь заменой условия $u|_{t=1}=f(x)$ на $u|_{t=0}=q(x)$. Граница $t=1$ --- наблюдаемая, $t=0$ --- ненаблюдаемая.

На пространстве непрерывно дифференцируемых функций со стандартным скалярным произведением $$H = C^2(\Pi)\cap L^2(\Pi)$$ определим оператор $A: H\rightarrow H$ следующим образом:

\begin{equation}
(Aq)(x)=u(x,1)\quad \forall x\in\Pi,
\label{heat3DA}
\end{equation}
где $u$ --- решение соответствующей корректной задачи Коши (\ref{heat3Du_correct}).

Вычисление сопряжённого оператора и градиента аналогично одномерному случаю, рассмотренному в статье \cite{pletnev2023component}, только вместо одного слагаемого с второй производной по $x$ --- три однотипных, для каждой пространственной переменной своё.

\begin{table}[H]
\begin{center}
\begin{tabular}{|c|c|c|}\hline
    $ $ & $||q_n-q^*||$ & $J(q_n)$ \\
   \hline
   ММО, $m=1$ & $1.65\times 10^{-3}$ & $1.70\times 10^{-14}$ \\
   \hline
   ММО, $m=2$ & $1.65\times 10^{-3}$ & $1.38\times 10^{-14}$ \\
   \hline
   ММО, $m=5$ & $1.64\times 10^{-3}$ & $7.27\times 10^{-15}$ \\
   \hline
   ММО, $m=\infty$ & $1.18\times 10^{-3}$ & $5.03\times 10^{-26}$ \\
   \hline
   Сопряжённые градиенты (ФР) & $1.66\times 10^{-3}$ & $7.51\times 10^{-16}$ \\
   \hline
   Сопряжённые градиенты (ПР) & $1.66\times 10^{-3}$ & $9.67\times 10^{-16}$ \\
   \hline
   Сопряжённые градиенты ($\perp$) & $6.59\times 10^{-3}$ & $1.46\times 10^{-9}$ \\
   \hline
   Тяжёлый шарик & $6.73\times 10^{-3}$ & $2.75\times 10^{-9}$ \\
   \hline
   Шаг Поляка & $3.09\times 10^{-3}$ & $1.76\times 10^{-12}$ \\
   \hline
   STM & $2.89\times 10^{-3}$ & $2.14\times 10^{-12}$ \\
   \hline
\end{tabular}
\end{center}
\caption{Результаты работы методов при решении решении задачи (\ref{heat3Du}) при $\kappa_{max}=0.4$ с использованием сетки с шагом $h=0.04$. Начальное значение функционала $J(q_0) = 6.27\times 10^{-3}$, начальное расстояние до точного решения $||q_0-q^*||=0.242$.}
\label{tab:table6.1}
\end{table}

\begin{figure}[H]
	\begin{minipage}[h]{0.5\linewidth}
		\center{\includegraphics[width=0.9\linewidth]{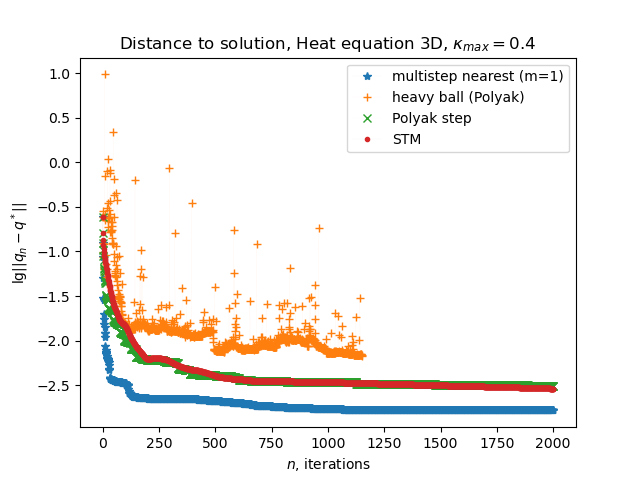}}
	\end{minipage}
	\hfill
	\begin{minipage}[h]{0.5\linewidth}
		\center{\includegraphics[width=0.9\linewidth]{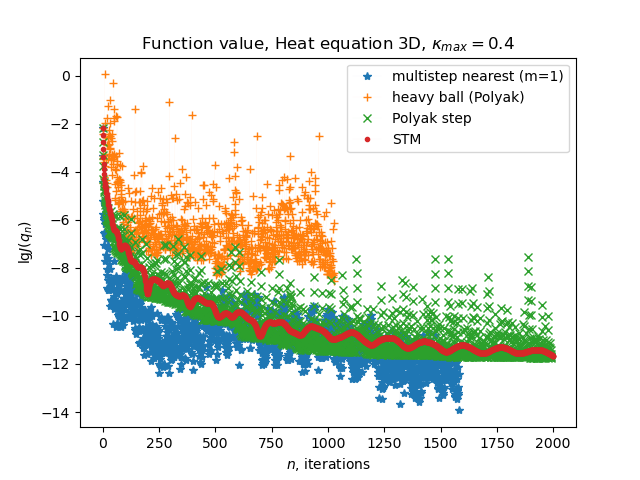}}
	\end{minipage}
	\caption{Сравнение $1$-моментного ММО, адаптивного тяжёлого шарика, градиентного спуска с шагом Поляка и метода подобных треугольников в применении к решению задачи (\ref{heat3Du}) при $\kappa_{max}=0.4$ с использованием сетки с шагом $h=0.04$.}
	\label{ris:image6.1}
\end{figure}

Рисунок \ref{ris:image6.1} показывает, что $1$-моментный метод минимальных ошибок позволяет достичь меньшей невязки, чем адаптивный метод тяжёлого шарика из статьи \cite{goujaud2022quadratic}, градиентный спуск с шагом Поляка и метод подобных треугольников, и по аргументу, и по функционалу. При этом сходимость по функционалу не является монотонной.

\begin{figure}[H]
	\begin{minipage}[h]{0.5\linewidth}
		\center{\includegraphics[width=0.9\linewidth]{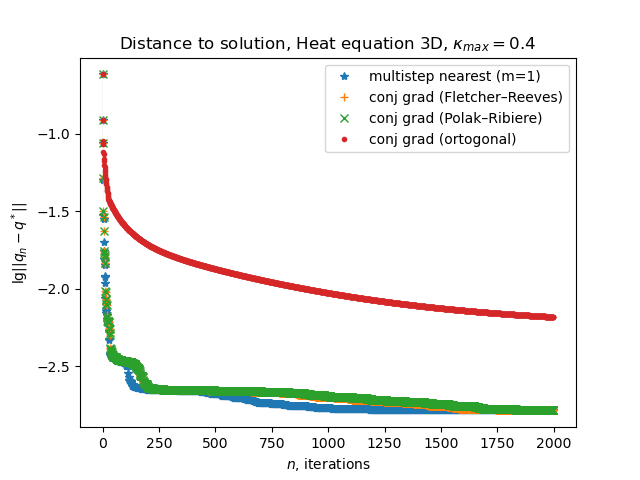}}
	\end{minipage}
	\hfill
	\begin{minipage}[h]{0.5\linewidth}
		\center{\includegraphics[width=0.9\linewidth]{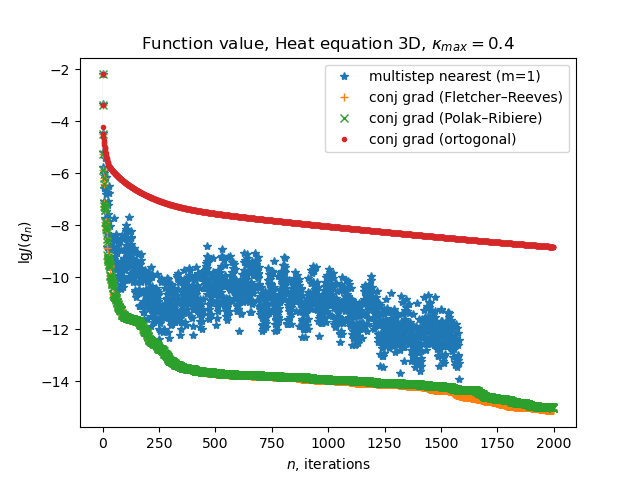}}
	\end{minipage}
	\caption{Сравнение $1$-моментного ММО и различных вариантов метода сопряжённых градиентов в применении к решению задачи (\ref{heat3Du}) при $\kappa_{max}=0.4$ с использованием сетки с шагом $h=0.04$.}
	\label{ris:image6.2}
\end{figure}

Рисунок \ref{ris:image6.2} показывает, что классические варианты метода сопряжённых градиентов и новый $1$-моментный метод минимальных ошибок дают близкую невязку по аргументу. Невязка по функционалу у метода сопряжённых градиентов оказалась меньше.

\begin{figure}[H]
	\begin{minipage}[h]{0.5\linewidth}
		\center{\includegraphics[width=1\linewidth]{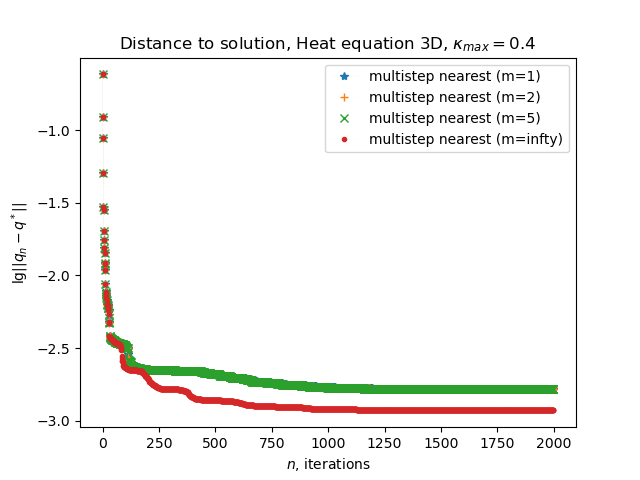}}
	\end{minipage}
	\hfill
	\begin{minipage}[h]{0.5\linewidth}
		\center{\includegraphics[width=1\linewidth]{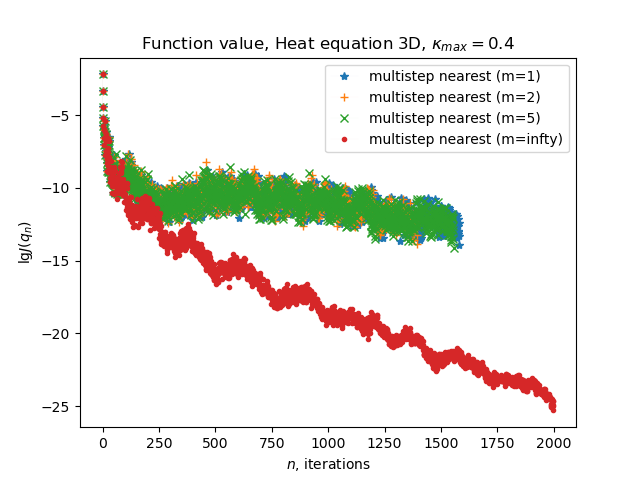}}
	\end{minipage}
	\caption{Сравнение $m$-моментных ММО с $m\in\{1, 2, 5, \infty\}$ в применении к решению задачи (\ref{heat3Du}) при $\kappa_{max}=0.4$ с использованием сетки с шагом $h=0.04$.}
	\label{ris:image6.3}
\end{figure}

Рисунок \ref{ris:image6.3} показывает, что $\infty$-моментный метод минимальных ошибок превосходит остальные рассмотренные методы: он достигает меньшей невязки и по аргументу, и по функционалу.

Эксперименты с $\kappa_{max}=0.6$ дают очень похожие результаты. Моментный метод минимальных ошибок достигает меньшего расстояния до точного решения, чем все остальные рассмотренные методы. $m=\infty$ позволяет достичь наименьшей невязки как по аргументу, так и по функционалу.

\begin{table}[H]
\begin{center}
\begin{tabular}{|c|c|c|}\hline
    $ $ & $||q_n-q^*||$ & $J(q_n)$ \\
   \hline
   ММО, $m=1$ & $3.74\times 10^{-3}$ & $9.93\times 10^{-16}$ \\
   \hline
   ММО, $m=2$ & $3.72\times 10^{-3}$ & $1.35\times 10^{-15}$ \\
   \hline
   ММО, $m=5$ & $3.72\times 10^{-3}$ & $2.36\times 10^{-15}$ \\
   \hline
   ММО, $m=\infty$ & $1.66\times 10^{-3}$ & $3.41\times 10^{-25}$ \\
   \hline
   Сопряжённые градиенты (ФР) & $4.76\times 10^{-3}$ & $3.89\times 10^{-16}$ \\
   \hline
   Сопряжённые градиенты (ПР) & $6.00\times 10^{-3}$ & $9.26\times 10^{-16}$ \\
   \hline
   Сопряжённые градиенты ($\perp$) & $74.3\times 10^{-3}$ & $4.14\times 10^{-8}$ \\
   \hline
   Тяжёлый шарик & $63.4\times 10^{-3}$ & $3.26\times 10^{-8}$ \\
   \hline
   Шаг Поляка & $26.0\times 10^{-3}$ & $1.19\times 10^{-10}$ \\
   \hline
   STM & $29.6\times 10^{-3}$ & $2.99\times 10^{-10}$ \\
   \hline
\end{tabular}
\end{center}
\caption{Результаты работы методов при решении решении задачи (\ref{heat3Du}) при $\kappa_{max}=0.6$ с использованием сетки с шагом $h=0.04$. Начальное значение функционала $J(q_0) = 1.34\times 10^{-3}$, начальное расстояние до точного решения $||q_0-q^*||=0.242$.}
\label{tab:table6.2}
\end{table}

\begin{figure}[H]
	\begin{minipage}[h]{0.5\linewidth}
		\center{\includegraphics[width=0.9\linewidth]{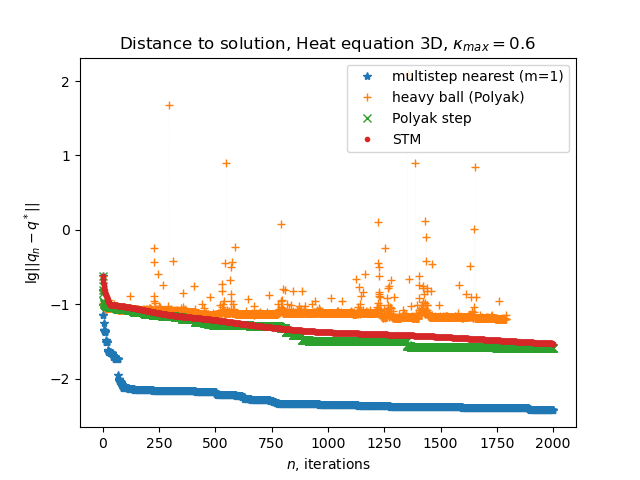}}
	\end{minipage}
	\hfill
	\begin{minipage}[h]{0.5\linewidth}
		\center{\includegraphics[width=0.9\linewidth]{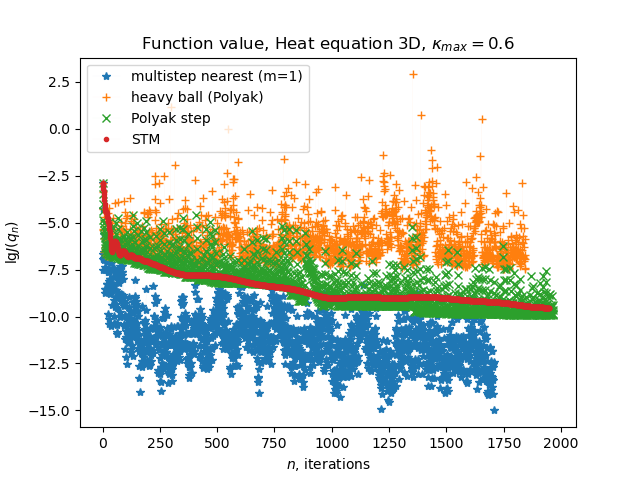}}
	\end{minipage}
	\caption{Сравнение $1$-моментного ММО, адаптивного тяжёлого шарика, градиентного спуска с шагом Поляка и метода подобных треугольников в применении к решению задачи (\ref{heat3Du}) при $\kappa_{max}=0.6$ с использованием сетки с шагом $h=0.04$.}
	\label{ris:image6.4}
\end{figure}

\begin{figure}[H]
	\begin{minipage}[h]{0.5\linewidth}
		\center{\includegraphics[width=0.9\linewidth]{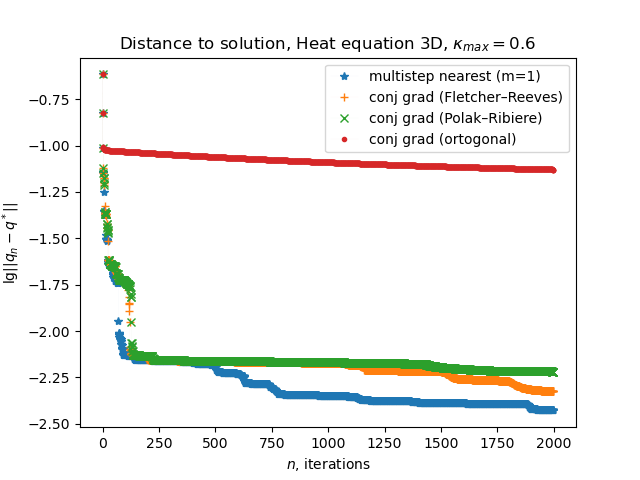}}
	\end{minipage}
	\hfill
	\begin{minipage}[h]{0.5\linewidth}
		\center{\includegraphics[width=0.9\linewidth]{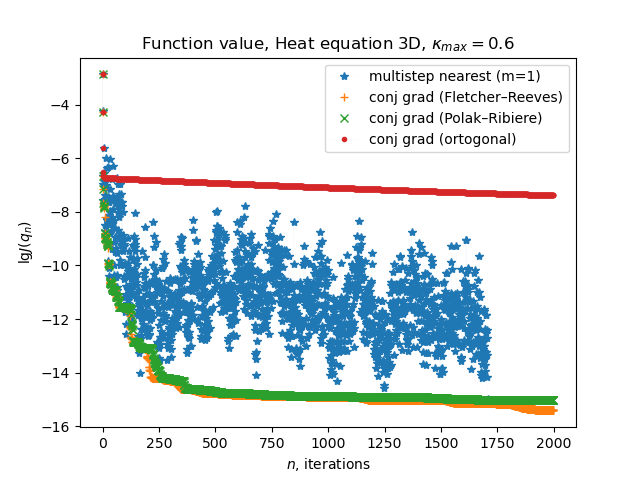}}
	\end{minipage}
	\caption{Сравнение $1$-моментного ММО и различных вариантов метода сопряжённых градиентов в применении к решению задачи (\ref{heat3Du}) при $\kappa_{max}=0.6$ с использованием сетки с шагом $h=0.04$.}
	\label{ris:image6.5}
\end{figure}

\begin{figure}[H]
	\begin{minipage}[h]{0.5\linewidth}
		\center{\includegraphics[width=0.9\linewidth]{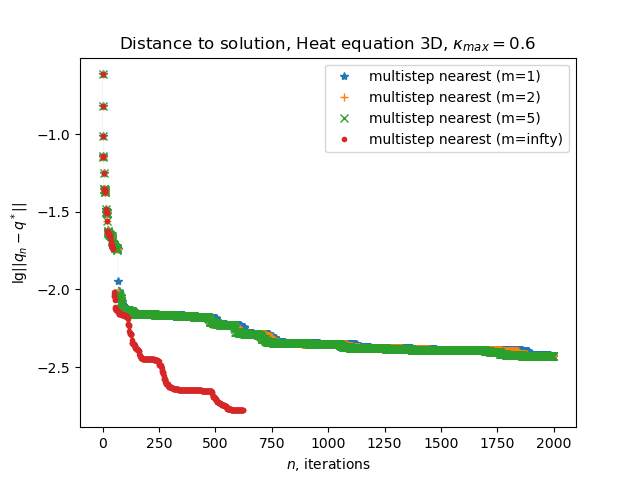}}
	\end{minipage}
	\hfill
	\begin{minipage}[h]{0.5\linewidth}
		\center{\includegraphics[width=0.9\linewidth]{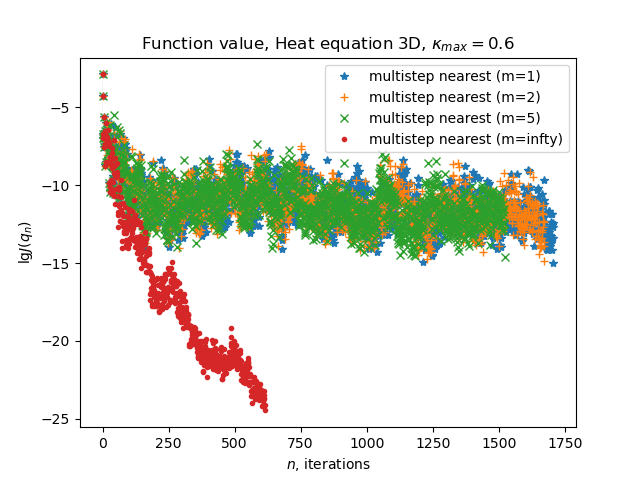}}
	\end{minipage}
	\caption{Сравнение $m$-моментных ММО с $m\in\{1, 2, 5, \infty\}$ в применении к решению задачи (\ref{heat3Du}) при $\kappa_{max}=0.6$ с использованием сетки с шагом $h=0.04$.}
	\label{ris:image6.6}
\end{figure}

Эксперименты показывают, что градиентные методы оптимизации позволяют эффективно решать ретроспективную задачу Коши для трёхмерного уравнения теплопроводности с граничными условиями на значение функции.

Увеличение коэффициента теплопроводности ухудшает качество решения. Это связано со свойством <<жёсткости>> спектра: собственные значения связанного с задачей самосопряжённого оператора имеют вид $\lambda_n=e^{-2\pi^2\kappa^2n^2}$, то есть стремятся к нулю значительно быстрее геометрической прогрессии. Чем больше $\kappa$, тем быстрее это стремление.

Моментный метод минимальных ошибок позволяет получить лучшее качество решения по сравнению с другими методами. Особенно хорошие результаты даёт метод с $m=\infty$.

Этот результат вполне согласуется с теоретическими выкладками --- теоремами 1 и  3.

\vspace{12pt}
\centerline{\it 5.3. Решение обратной задачи термоакустики}
\vspace{12pt}

В статье \cite{kabanikhin2011krivorotko} исследуется обратная задача термоакустики. Рассмотрим обратную задачу $3$ с $L=T=1$:

\begin{equation}
\begin{cases}u_{tt} = u_{xx} + u_{yy}, & (x,y)\in \Omega = (0,1)\times (0,1), t\in (0, 1)\\u_t|_{t=0} = 0, & (x,y)\in \Omega\\
    u_x|_{x=0} = u_x|_{x=1} = 0, & y\in (0,1), t\in (0, 1)\\
    u_y|_{y=0} = u_y|_{y=1} = 0, & x\in (0,1), t\in (0, 1)\\
    u|_{x=0}=f_1(y,t), & y\in (0,1), t\in (0, 1)\\
    u|_{x=1}=f_2(y,t), & y\in (0,1), t\in (0, 1)\\
    u|_{y=1}=f_3(x,t), & x\in (0,1), t\in (0, 1)\end{cases}
\label{3thermoacoustic}
\end{equation}

Как доказано в статье \cite{kabanikhin2011krivorotko}, эта задача некорректна. Корректной является прямая задача термоакустики
\begin{equation}
\begin{cases}u_{tt} = u_{xx} + u_{yy}, & (x,y)\in \Omega = (0,1)\times (0,1), t\in (0, 1)
\\u|_{t=0} = q(x,y), & (x,y)\in \Omega
\\u_t|_{t=0} = 0, & (x,y)\in \Omega\\
    u_x|_{x=0} = u_x|_{x=1} = 0, & y\in (0,1), t\in (0, 1)\\
    u_y|_{y=0} = u_y|_{y=1} = 0, & x\in (0,1), t\in (0, 1)\end{cases}
\label{thermoacoustic_correct},
\end{equation}

которая получается заменой граничных условий на функцию начальным условием. По смыслу задачи, именно начальное условие нужно найтя, зная граничные. Функция $u(x,y,t)$ определена на параллелепипеде. Его граница $t=0$ --- ненаблюдаемая; границы $x=0$, $x=1$, $y=1$ --- наблюдаемые.

Для каждой наблюдаемой границы можно определить свой оператор. Важно, что все они действуют из $C^1([0, L]\times [0, L])$ в $C^1([0, L]\times [0, T])$. В общем случае, области определения функций $q$ и $A_iq$ могут оказаться различными.

\numberwithin{equation}{subsection}
\begin{equation}
\begin{gathered}
(A_1q)(y, t) = u(0, y, t)\quad \forall (y, t)\in [0, 1]\times [0, 1]
\\
(A_2q)(y, t) = u(1, y, t)\quad \forall (y, t)\in [0, 1]\times [0, 1]
\\
(A_3q)(x, t) = u(x, 1, t)\quad \forall (x, t)\in [0, 1]\times [0, 1]
\label{thermoacoustic_A}
\end{gathered}
\end{equation}

Во всех случаях $u(x, y, t)$ --- решение корректной задачи (\ref{thermoacoustic_correct}).

Функционал является составным:
\begin{equation}
J(q) = \frac{1}{2}\sum\limits_{l=1}^3 ||A_lq-f_l||^2
\label{thermoacoustic_J}
\end{equation}

Градиент, как обычно, вычисляется с помощью сопряжённых операторов. Приведём их формулы согласно \cite{kabanikhin2011krivorotko}:
$$A_l: C^1([0, 1]\times [0,1])\rightarrow C^1([0, 1]\times [0,1])$$
$$(A_1^*p)(x, y)=\psi_{1t}(x, y, 0), \text{ где } \psi_1(x, y, t)\text{ - решение задачи}$$
$$\begin{cases}\psi_{1tt} = \psi_{1xx} + \psi_{1yy}, & (x,y)\in \Omega = (0,1)\times (0,1), t\in (0, 1)
\\\psi_1|_{t=T} = \psi_{1t}|_{t=T} = 0, & (x,y)\in \Omega\\
    \psi_{1x}|_{x=0} = p(y, t), & y\in (0,1), t\in (0, 1)\\
    \psi_{1x}|_{x=1} = 0, & y\in (0,1), t\in (0, 1)\\
    \psi_{1y}|_{y=0} = \psi_{1y}|_{y=1} = 0, & x\in (0,1), t\in (0, 1)\end{cases}$$

$$(A_2^*p)(x, y)=\psi_{2t}(x, y, 0), \text{ где } \psi_2(x, y, t)\text{ - решение задачи}$$
$$\begin{cases}\psi_{2tt} = \psi_{2xx} + \psi_{2yy}, & (x,y)\in \Omega = (0,1)\times (0,1), t\in (0, 1)
\\\psi_2|_{t=T} = \psi_{2t}|_{t=T} = 0, & (x,y)\in \Omega\\
    \psi_{2x}|_{x=0} = 0, & y\in (0,1), t\in (0, 1)\\
    \psi_{2x}|_{x=1} = -p(y, t), & y\in (0,1), t\in (0, 1)\\
    \psi_{2y}|_{y=0} = \psi_{2y}|_{y=1} = 0, & x\in (0,1), t\in (0, 1)\end{cases}$$

$$(A_3^*p)(x, y)=\psi_{3t}(x, y, 0), \text{ где } \psi_2(x, y, t)\text{ - решение задачи}$$
$$\begin{cases}\psi_{3tt} = \psi_{3xx} + \psi_{3yy}, & (x,y)\in \Omega = (0,1)\times (0,1), t\in (0, 1)
\\\psi_3|_{t=T} = \psi_{3t}|_{t=T} = 0, & (x,y)\in \Omega\\
    \psi_{3x}|_{x=0} = \psi_{3x}|_{x=1} = 0, & y\in (0,1), t\in (0, 1)\\
    \psi_{3y}|_{y=0} = 0, & x\in (0,1), t\in (0, 1)\\
    \psi_{3y}|_{y=1} = -p(x, t), & x\in (0,1), t\in (0, 1)\end{cases}$$

$$\nabla J(q) = \sum\limits_{l=1}^3 A_l^*\left(A_lq-f_l\right).$$

Для проверки методов используется задача с $$q^*(x, y) = 0.1+\frac{(1+\cos{8\pi x})(1+\cos{8\pi y})}{32}\left[x, y\in \left[\frac{1}{8}, \frac{3}{8}\right]\cup \left[\frac{5}{8}, \frac{7}{8}\right]\right].$$
Внешние квадратные скобки, как обычно, обозначают функцию-индикатор, которая принимает значение $1$, если утверждение внутри скобок истинно, и $0$ в противоположном случае. На такой же задаче проводились эксперименты в статье \cite{kabanikhin2011krivorotko}.

$q^*(x,y)$ моделирует область с четырьмя неоднородностями.

Для приближённого решения задач применяются сетки с $h=0.02$, $\tau=0.002$.

\begin{table}[H]
\begin{center}
\begin{tabular}{|c|c|c|}\hline
    $ $ & $||q_n-q^*||$ & $J(q_n)$ \\
   \hline
   ММО, $m=1$ & $8.35\times 10^{-5}$ & $6.93\times 10^{-12}$ \\
   \hline
   ММО, $m=2$ & $5.98\times 10^{-5}$ & $9.47\times 10^{-13}$ \\
   \hline
   ММО, $m=5$ & $4.07\times 10^{-5}$ & $2.29\times 10^{-13}$ \\
   \hline
   ММО, $m=\infty$ & $3.32\times 10^{-3}$ & $8.65\times 10^{-6}$ \\
   \hline
   Сопряжённые градиенты (ФР) & $3.03\times 10^{-3}$ & $3.90\times 10^{-6}$ \\
   \hline
   Сопряжённые градиенты (ПР) & $8.42\times 10^{-5}$ & $1.47\times 10^{-12}$ \\
   \hline
   Сопряжённые градиенты ($\perp$) & $1.96\times 10^{-4}$ & $3.63\times 10^{-11}$ \\
   \hline
   Тяжёлый шарик & $3.74\times 10^{-3}$ & $1.18\times 10^{-5}$ \\
   \hline
   Шаг Поляка & $5.59\times 10^{-5}$ & $2.77\times 10^{-13}$ \\
   \hline
   STM & $7.33\times 10^{-3}$ & $3.91\times 10^{-5}$ \\
   \hline
\end{tabular}
\end{center}
\caption{Результаты работы методов при решении решении задачи (\ref{3thermoacoustic}) с использованием сетки с шагами $h=0.02$, $\tau=0.002$. Начальное значение функционала $J(q_0) = 0.018$, начальное расстояние до точного решения $||q_0-q^*||=0.11$.}
\label{tab:table7.3}
\end{table}

\begin{figure}[H]
	\begin{minipage}[h]{0.5\linewidth}
		\center{\includegraphics[width=0.9\linewidth]{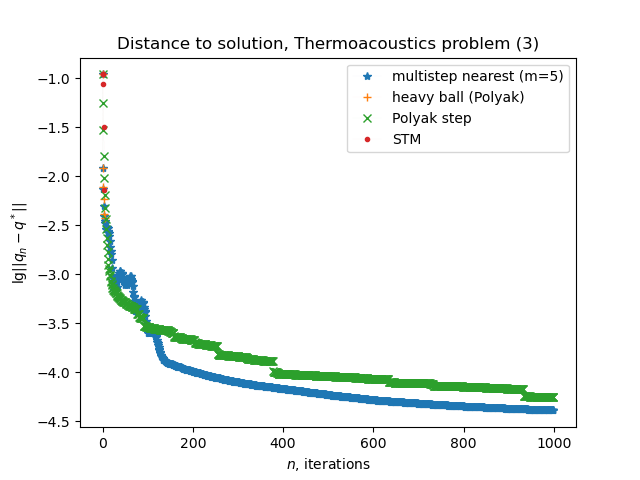}}
	\end{minipage}
	\hfill
	\begin{minipage}[h]{0.5\linewidth}
		\center{\includegraphics[width=0.9\linewidth]{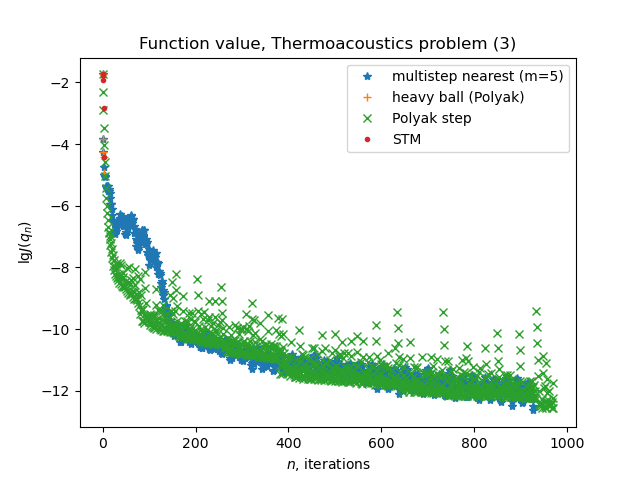}}
	\end{minipage}
	\caption{Сравнение $5$-моментного ММО, адаптивного тяжёлого шарика, градиентного спуска с шагом Поляка и метода подобных треугольников в применении к решению задачи (\ref{3thermoacoustic}) с использованием сетки с шагами $h=0.02$, $\tau=0.002$.}
	\label{ris:image7.7}
\end{figure}

\begin{figure}[H]
	\begin{minipage}[h]{0.5\linewidth}
		\center{\includegraphics[width=0.9\linewidth]{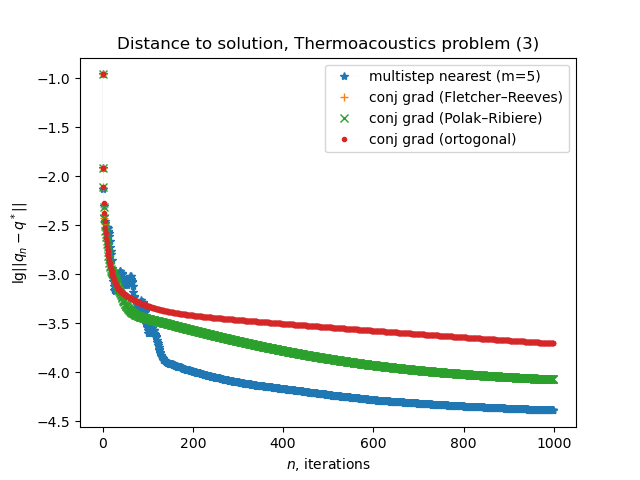}}
	\end{minipage}
	\hfill
	\begin{minipage}[h]{0.5\linewidth}
		\center{\includegraphics[width=0.9\linewidth]{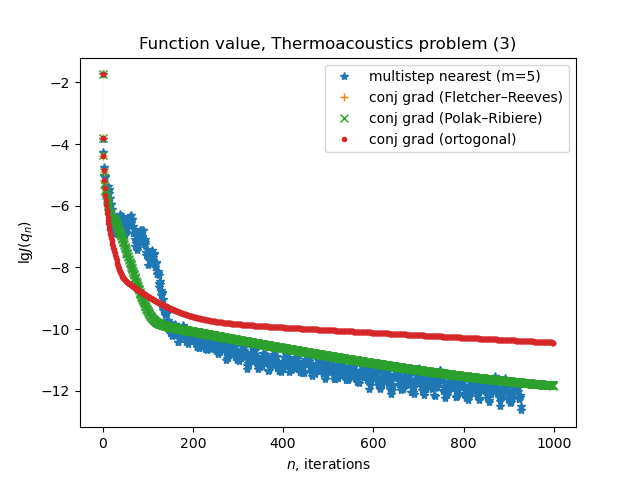}}
	\end{minipage}
	\caption{Сравнение $5$-моментного ММО и различных вариантов метода сопряжённых градиентов в применении к решению задачи (\ref{3thermoacoustic}) с использованием сетки с шагами $h=0.02$, $\tau=0.002$.}
	\label{ris:image7.8}
\end{figure}

Рисунок \ref{ris:image7.7} показывает, что только $5$-моментный метод минимальных ошибок и градиентный спуск с шагом Поляка демонстрируют монотонную сходимость. При этом $5$-моментный метод минимальных ошибок показывает значительно лучшие результаты.

Метод сопряжённых градиентов в форме Флетчера-Ривса перестаёт уменьшать невязку после нескольких итераций. Метод сопряжённых градиентов в форме Полака-Рибьера, как и метод с ортогональными шагами и вспомогательной минимизацией функционала, продолжает работать, но его результаты хуже, чем у $5$-моментного ММО.

\begin{figure}[H]
	\begin{minipage}[h]{0.5\linewidth}
		\center{\includegraphics[width=0.9\linewidth]{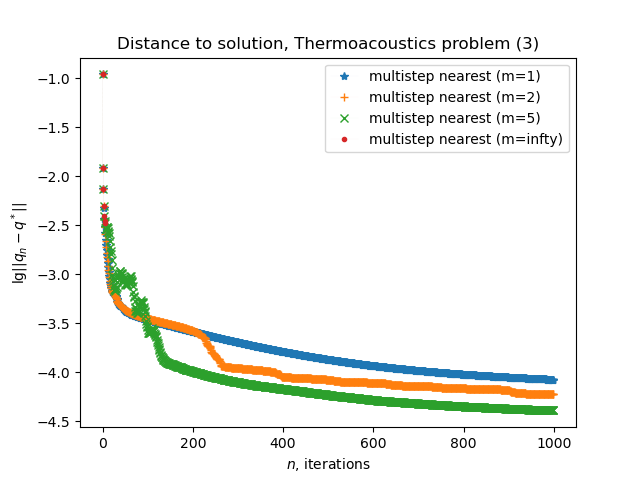}}
	\end{minipage}
	\hfill
	\begin{minipage}[h]{0.5\linewidth}
		\center{\includegraphics[width=0.9\linewidth]{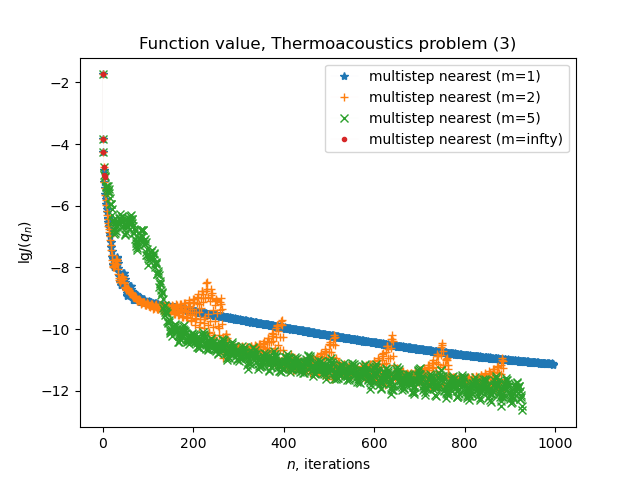}}
	\end{minipage}
	\caption{Сравнение $m$-моментных ММО с $m\in\{1, 2, 5, \infty\}$ в применении к решению задачи (\ref{3thermoacoustic}) с использованием сетки с шагами $h=0.02$, $\tau=0.002$.}
	\label{ris:image7.9}
\end{figure}

Среди $m$-моментных методов минимальных ошибок наилучшие результаты достигнуты при $m=5$. Метод с $m=\infty$, наилучший в теории, слишком чувствителен к погрешностям вычислений.

$5$-моментный метод минимальных ошибок во всех случаях показал наилучшие результаты. Это подтверждает выводы раздела 3.4 о возможности применения построенных методов к составному функционалу --- сумме квадратичных функционалов, имеющих общую точку минимума.

Теоретические ожидания наилучшей сходимости при $m=\infty$ не оправдались. Это объясняется чувствительностью такого метода к погрешностям и их накоплением.

Эксперименты показывают возможность решения обратной задачи термоакустики с помощью градиентных методов оптимизации, в особенности --- нового $m$-моментного метода минимальных ошибок.

\vspace{12pt}
\centerline{6. Заключение} 
\vspace{12pt}

В статье построен $m$-моментный метод минимальных ошибок для решения задач оптимизации квадратичных функционалов. Доказана его оптимальность (теорема 1): никакой метод первого порядка не может сходиться по аргументу быстрее. Также установлен факт сходимости даже в бесконечномерном пространстве (теорема 3). При этом сходимость не является равномерной, о чём свидетельствует теорема 6.

Численные эксперименты показывают высокую эффективность применения построенного метода к решению различных некорректных задач для дифференциальных уравнений, которые сводятся именно к квадратичным задачам оптимизации в гильбертовом пространстве.

\end{document}